\documentclass[12pt]{article}
\usepackage{graphics}
\newtheorem{THEOREM}{Theorem}[section]
\newtheorem{LEMMA}{Lemma}[section]
\newtheorem{COROLLARY}{Corollary}[section]

\newcommand{\qed}{\ \rule[-1pt]{4pt}{8pt} 

                                         \vspace{2ex} }
\newenvironment{PROOF}{

                       \noindent{\bf Proof}.}{\qed}
\newcounter{labelflag} 

\newcommand{\Label}[1]{
                       \ifnum\thelabelflag=1 
                          \ifmmode  
                             \makebox[0in][l]{\qquad\fbox{\rm#1}}
                          \else
                             \marginpar{\vspace{0.7\baselineskip}
                                        \hspace{-1.1\textwidth}
                                        \fbox{\rm#1}}
                          \fi 
                       \fi
                       \label{#1} 
                      }
\newcommand{\be}{\begin{equation}}
\newcommand{\ee}{\end{equation}}
\newcommand{\eps}{\varepsilon}
\newcommand{\half}{{\textstyle{\frac{1}{2}}}}

\newcommand{\onefourth}{{\textstyle{\frac{1}{4}}}}
\newcommand{\onesixth}{{\textstyle{\frac{1}{6}}}}
\begin{document}

\vspace*{0.0in}
\begin{center}
{\large\textbf{ASYMPTOTIC ANALYSIS OF}} \\[1ex]
{\large\textbf{TWO REDUCTION METHODS FOR}} \\[1ex]
{\large\textbf{SYSTEMS OF CHEMICAL REACTIONS}}

\vspace{2ex}
Hans G.\ Kaper$^1$ and
Tasso J.~Kaper$^2$

$^1$~Mathematics and Computer Science Division \\
Argonne National Laboratory,
Argonne, Illinois, USA \\[1ex]

$^2$~Department of Mathematics and Center for BioDynamics\\
Boston University,
Boston, Massachusetts, USA 
\end{center}

\begin{small}
\paragraph{Abstract.}
This article concerns
two methods for reducing
large systems of chemical
kinetics equations,
namely, the method of
intrinsic low-dimensional
manifolds (ILDMs) due to
Maas and Pope
[U.~Maas and S.~B.~Pope, \emph{Combustion and Flame} \textbf{88} (1992) 239--264]
and an iterative method due to
Fraser
[S.~J.~Fraser, \emph{J.~Chem.\ Phys.} \textbf{88} (1988) 4732--4738]
and further developed by
Roussel and Fraser
[M.~R.~Roussel and S.~J.~Fraser, \emph{J.~Chem.\ Phys.} \textbf{93} (1990) 1072--1081].
Both methods exploit
the separation of fast
and slow reaction time scales
to find low-dimensional manifolds
in the space of species concentrations
where the long-term dynamics are played out.
The analysis is carried out 
in the context of systems of
ordinary differential equations
with multiple time scales and
geometric singular perturbation
theory (GSPT).
A small parameter $\eps$
measures the separation 
of time scales.
The underlying assumption is
that the system of equations
has an asymptotically stable 
slow manifold $\mathcal{M}_0$
in the limit as $\eps \downarrow 0$.
Then it follows from GSPT
that there exists a slow manifold
$\mathcal{M}_\eps$ for all
sufficiently small positive $\eps$,
which is asymptotically close
to~$\mathcal{M}_0$.

It is shown that the ILDM method
yields a low-dimensional manifold
whose asymptotic expansion agrees
with the asymptotic expansion
of $\mathcal{M}_\eps$
up to and including terms
of $\mathcal{O} (\eps)$.
At $\mathcal{O} (\eps^2)$,
an error appears that is
proportional to
the local curvature
of $\mathcal{M}_0$;
it vanishes if and only if
the curvature is zero
everywhere.

The iterative method
generates, term by term,
the asymptotic expansion
of the slow manifold~$\mathcal{M}_\eps$.
Starting from $\mathcal{M}_0$,
the $i$th application
of the algorithm yields
the correct expansion coefficient
at $\mathcal{O}(\eps^i)$,
while leaving the lower-order
coefficients at $\mathcal{O}(1)$
through $\mathcal{O} (\eps^{i-1})$
invariant.
Thus, after $\ell$ applications,
the expansion is accurate
up to and including the terms
of $\mathcal{O}(\eps^\ell)$.

The analytical results are
illustrated in two examples:
a planar system from
enzyme kinetics
(Michaelis--Menten--Henri)
and a model planar system
due to Davis and Skodje.

\noindent\textbf{AMS Classification (MSC2000).}
Primary: 80A30, 34E13.
Secondary: 92C45, 34C20, 34E15.



\end{small}
%
\section{Introduction and Summary of Results \label{s-intro}}
\setcounter{equation}{0}
Many chemical reaction mechanisms in
combustion~\cite{W-1985, PR-1993, MDWZ-1999},
atmospheric science~\cite{SP-1998},
enzyme kinetics~\cite{F-1975},
and biochemistry~\cite{KF-2001}
involve large numbers of species,
multiple chains of chemical reactions,
and widely disparate time scales.
A typical model of hydrocarbon combustion,
for example, may well involve several hundred
species, which participate in hundreds of
reactions that proceed on time scales
ranging from nanoseconds to minutes.
The size and complexity of these mechanisms
has stimulated the search for methods
that reduce the number of species and
chemical reactions but retain a desired
degree of accuracy.
Typically, these \emph{reduction methods}
select a small number of species,
which are marked as
\emph{reaction progress variables},
and determine the concentrations
of the remaining species
as functions of the latter,
either by table look-ups or
by direct computation.
The critical step in
these methods is, of course,
the definition of the reaction
progress variables, which may be
actual concentrations of
selected species or
combinations thereof.

Research into reduction methods
has increased dramatically 
over the past decade, and
several methods have been proposed
in the literature and implemented
in computer codes.
We mention
the quasi-steady state
approximation~\cite{SS-1989, TTP-1993, YTBMP-1996},
the partial-equilibrium approximation~\cite{S-1991},
methods based on details of the chemistry~\cite{S-1991, PR-1993},
an iterative method~\cite{F-1988, RF-1990},
the method of intrinsic low-dimensional
manifolds~\cite{MP-1992conf, MP-1992},
the computational singular perturbation
method~\cite{GL-1992, L-1993, LG-1994, HG-1999, MDMG-1999, LJL-2001},
a principal-component analysis,
lumping techniques~\cite{LTRT-1994},
repro-modeling~\cite{T-1994},
an inertial-manifold approach~\cite{YTBMP-1995},
a dynamic dimension-reduction method~\cite{DH-1996, DHM-1996},
a saddle-point method~\cite{DS-1999, DS-2001},
a predictor-corrector method~\cite{DS-1999, DS-2001},
an optimization method~\cite{PZ-1999},
and a global-eigenvalue method~\cite{SD-2001}.

In this article, we focus on two
reduction methods, namely, the
intrinsic low-dimensional manifold
(ILDM) method due to Maas and
Pope~\cite{MP-1992conf, MP-1992}
and the iterative method due to Fraser
and further developed by
Roussel and Fraser~\cite{F-1988, RF-1990}.
Both methods have been developed for
and extensively applied to problems
with slow manifolds that attract
nearby initial conditions.
The long-time behavior of
such systems is governed by
the dynamics on the slow manifold,
whose dimension is generally 
much less than that of the
total composition space,
resulting in a considerable
reduction of complexity.

Given the importance of slow
manifolds, a central question
for any reduction method is:
How accurately does it
approximate a slow manifold?
The present investigation
answers this question
for the ILDM method
of Maas and Pope
and the iterative method
of Fraser and Roussel.

In the ILDM method,
the Jacobian of the vector field
is partitioned at each point of
phase space into a fast and a
slow component, and bases
for the corresponding subspaces
are generated by means
of a Schur decomposition.
The intrinsic low-dimensional
manifold is defined as
the locus of points where
the vector field lies entirely
in the slow subspace and is an
approximation of the slow manifold.
The efficacy of the ILDM method
is evident, for example,
by the reduction achieved
in the prototypical example
of a CO-H$_2$-O$_2$-N$_2$
combustion model~\cite{MP-1992, M-1998}.
Disregarding only the production of NO,
the model comprises evolution equations
for the enthalpy,
pressure, and concentration of
each of thirteen species,
making for a fifteen-dimensional
phase space,
and a total of sixty-seven chemical reactions.
With the proper choice of the reaction
progress variable (CO$_2$),
a reduction to a one-dimensional
ILDM can be achieved that retains
a certain accuracy after an
initial transient~\cite{MP-1992, M-1998}.
Reduction to a two-dimensional ILDM
gives a better approximation,
albeit at the expense
of keeping track of 
two reaction progress variables
and the storage of a correspondingly
larger look-up table.
Refinements, applications,
and evaluations of the ILDM method
against direct numerical simulations
can be found in Refs.~\cite{MP-1994,
M-1995, ED-1995, SM-1996, SMSRW-1996,
ELKD-1997, BSM-1997, NSM-1997, YP-1998a, YP-1998b,
BSM-1998, SBM-1998, RMW-1999,
BBM-2000, BM-2000, YM-2000,
M-2001, SRPP-2001}.

The iterative method was
inspired by the phase-space geometry
of an enzyme-kinetics model involving
a fast and a slow species,
where the slow manifold
is a curve in the phase plane.
The method is derived formally from
the \emph{invariance equation}---an
equation that is satisfied
on any trajectory of
the dynamical system and,
in particular, on the slow manifold
and extends naturally to
multidimensional systems with
(possibly) higher-dimensional
slow manifolds.
The procedure is explicit
if the force field is linear
in the fast variable,
and implicit otherwise;
hence, it generally requires
the use of a nonlinear equation solver.
The method has been developed
further and applied,
in particular, to several problems
of enzyme kinetics and
metabolism in Refs.~\cite{NF-1989,
R-1994, R-1997, F-1998,
RF-1991a, RF-1991b, RF-1993, RF-2001}.

A natural framework for the analysis
of these and similar reduction methods
is provided by geometric singular
perturbation theory (GSPT)
\cite{F-1979, S-1990, J-1994, K-1999}.
The presence of a fast and a slow
time scale leads naturally to
the introduction of a small
positive parameter~$\eps$
measuring the ratio of
the characteristic times.
If, in the limit as $\eps \downarrow 0$
(infinite separation of time scales),
the system of kinetics equations
has a slow manifold, $\mathcal{M}_0$,
in phase space and this manifold is
asymptotically stable, then GSPT
identifies a (usually nonunique)
slow manifold $\mathcal{M}_\eps$
for $\eps$ sufficiently small positive
and gives a complete geometric and
analytical description of
all solutions in the vicinity
of the slow manifold,
including how trajectories
approach the manifold.
By comparing the asymptotics of
the slow manifold $\mathcal{M}_\eps$
found by GSPT
with the asymptotics of the
low-dimensional manifolds generated by
the ILDM method and the iterative method
we can evaluate the accuracy of
these reduction methods
for small values of $\eps$
(finite but large separation
of time scales).
The evaluation leads to
the following conclusions.
\begin{description}
\item [ILDM Method.]
(i)~The asymptotic expansion
of the ILDM agrees with the
asymptotic expansion of the
slow manifold $\mathcal{M}_\eps$
up to and including 
the $\mathcal{O}(\eps)$ term,
for all fast-slow systems.
In general, however,
the $\mathcal{O}(\eps^2)$ terms differ. \\[1ex]
(ii)~The error at $\mathcal{O} (\eps^2)$
is proportional to the local curvature 
of the slow manifold~$\mathcal{M}_0$.
It vanishes if and only if the curvature
of $\mathcal{M}_0$ is zero everywhere.
(The ``if'' part was observed
previously in Ref.~\cite{MP-1992}.)
\item [Iterative Method.]
(i)~The iterative method, if
started from~$\mathcal{M}_0$,
generates term by term
the asymptotic expansion 
of the slow manifold~$\mathcal{M}_\eps$.
In particular,
$\ell$ applications of the
iterative method generate
an approximation to the
slow manifold $\mathcal{M}_\eps$
that is asymptotically correct
up to and including the
$\mathcal{O}(\eps^\ell)$ term,
albeit with extraneous terms at
$\mathcal{O}(\eps^{\ell+1})$. \\[1ex]
(ii)~The $\ell$th iteration leaves
the terms at $\mathcal{O}(1)$
through $\mathcal{O}(\eps^{\ell-1})$
invariant.
(This observation is important
because the lower-order terms
have already been determined correctly
in the preceding iterations.)
\end{description}
\noindent\textbf{Remark.}
In Ref.~\cite{RMW-1999},
it is shown that the ILDM coincides with
the slow manifold $\mathcal{M}_0$
in the limit of infinite separation 
of the fast and slow time scales
($\eps=0$).~\qed

\noindent\textbf{Remark.}
The slow manifold $\mathcal{M}_0$
can often be found analytically;
otherwise, it can be obtained by one
application of the iterative method
to the steady-state approximation.
The latter is readily found numerically;
see Refs.~\cite{F-1988,RF-1990}.~\qed

The present article
is organized as follows.
In Section~\ref{s-general},
we review the general framework
of fast-slow systems of ordinary
differential equations and
recall the asymptotic expansion
of the slow manifold.
In Section~\ref{s-ILDM},
we define 
the ILDM 
and indicate briefly how it is computed.
We present the asymptotic expansion of the ILDM
for planar fast-slow systems
(one fast and one slow variable)
in Section~\ref{s-ILDM-planar}
and for general fast-slow systems
($n$ fast and $m$ slow variables)
in Section~\ref{s-ILDM-general}.
The results are summarized in
Corollary~\ref{c-ILDM}.
In Section~\ref{s-FR},
we describe the iterative method
of Fraser and Roussel.
We discuss its asymptotics
in Section~\ref{s-FR-as}.
The results are summarized
in Corollary~\ref{c-FR}.
We illustrate
the analytical results
with two planar examples, namely
the Michaelis--Menten--Henri mechanism
of enzyme kinetics (Section~\ref{s-MMH})
and an example due to Davis and Skodje
(Section~\ref{s-DS}).
In Section~\ref{s-disc}
we remark on several generalizations
and discuss some remaining issues.

\section{Fast-Slow Systems of ODEs \label{s-general}}
\setcounter{equation}{0}
We consider reaction mechanisms
in homogeneous media,
where the concentrations
of the chemical species
depend on time only.
The concentrations evolve
on two distinct and widely
separated time scales.
The slowly evolving concentrations are
the entries of the vector~$y$,
the remaining concentrations
the entries of the vector~$z$;
the former has $m$ components,
the latter~$n$ ($m, n \geq 1$).
The separation of time scales
is measured by $\eps$,
an arbitrarily small
positive parameter.
The limit $\eps \downarrow 0$
corresponds to infinite separation.
The reaction mechanism is thus modeled
by a system of ordinary differential
equations (ODEs),
\begin{eqnarray}
  y' &=& \eps f (y, z, \eps) ,  \Label{eq-y} \\
  z' &=& g (y, z, \eps) . \Label{eq-z}
\end{eqnarray}
The unknowns $y$ and $z$ are functions of $t$
with values in $\mathbf{R}^m$ and $\mathbf{R}^n$,
respectively;
${}'$ denotes differentiation with respect to $t$;
and $f$ and $g$ are smooth functions
with values in $\mathbf{R}^m$ and
$\mathbf{R}^n$, respectively.
We assume that $f$ and $g$,
as well as all their derivatives,
are $\mathcal{O} (1)$
as $\eps \downarrow 0$.

\noindent\textbf{Remark.}
The system of Eqs.~(\ref{eq-y})--(\ref{eq-z})
is, of course, an idealization
of the complex systems
that occur in chemical kinetics.
The model is adopted here
because it is suitable
for mathematical analysis.
We claim, however, that it also captures
the essential elements of any
reaction mechanism whose long-term
dynamics evolve on slow manifolds
and offers a paradigm for the analysis
of reduction methods.
The validity of our conclusions
extends therefore well beyond the
idealized system of
Eqs.~(\ref{eq-y})--(\ref{eq-z}).
We comment on the implications
for more realistic systems
in Section~\ref{s-disc}.~\qed

The independent variable~$t$
is called the \emph{fast time}
because it defines the time scale 
on which the fast variables evolve,
and the system of Eqs.~(\ref{eq-y})--(\ref{eq-z})
is labeled the \emph{fast system}.
While the fast time scale is
appropriate for the study of
the transient dynamics,
the long-time dynamics
are more naturally studied
in terms of the
\emph{slow time}
$\tau=\eps t$.
On the scale of $\tau$,
the system of
Eqs.~(\ref{eq-y})--(\ref{eq-z}) 
assumes the form
\begin{eqnarray}
  \dot{y} &=& f (y, z, \eps) ,
  \Label{eq-y-slow} \\
  \eps \dot{z} &=& g (y, z, \eps) .
  \Label{eq-z-slow}
\end{eqnarray}
Here, $\dot{~}$ denotes 
differentiation 
with respect to $\tau$.
We refer to the system of
Eqs.~(\ref{eq-y-slow})--(\ref{eq-z-slow})
as the \emph{slow system}.

The fast system~(\ref{eq-y})--(\ref{eq-z})
and the slow system~(\ref{eq-y-slow})--(\ref{eq-z-slow})
are, of course, equivalent 
as long as $\eps > 0$,
but they approach different limits
as $\eps \downarrow 0$---that is,
as the separation
of the fast and slow time scales
becomes infinite.
The fast system reduces to
\begin{eqnarray}
  y' &=& 0 ,  \Label{fast-red-y} \\
  z' &=& g (y, z, 0)  \Label{fast-red-z} ,
\end{eqnarray}
which is essentially a single equation 
for the fast variable $z$ with $y$ as a parameter.
The slow system, on the other hand,
reduces to
\begin{eqnarray}
  \dot{y} &=& f(y, z, 0) , \Label{slow-red-y} \\
  0 &=& g (y, z, 0) \Label{slow-red-z} .
\end{eqnarray}
The first equation describes
the motion of the slow variable $y$,
and the second equation is
an algebraic constraint
that forces the motion to take place
on the zero set of $g$.

Our focus is on systems for which
the zero set of $g$ is represented
by the graph of a function.
That is, we assume that there exists
a single-valued function $h_0$,
which is defined on a compact domain
$K = [0,Y]^m$ in $\mathbf{R}^m$,
such that
\be
  g (y, h_0(y), 0) = 0 , \quad y \in K .
  \Label{h0}
\ee
The zero set of $g$ thus defines
a \emph{manifold}, $\mathcal{M}_0$,
in phase space,
\be
  \mathcal{M}_0 = \{ (y, z) \in \mathbf{R}^{m+n} :
  z = h_0 (y), \; y \in K \} ,
  \Label{M0}
\ee
to which the motion of the reduced
slow system is confined.

Our analysis requires a second assumption
that holds for many, though not all,
of the systems
in which reductions have been sought,
namely, that each point $(y, h_0(y))$
on $\mathcal{M}_0$
is an \emph{asymptotically stable fixed point}
of Eq.~(\ref{fast-red-z}).
The assumption guarantees that
the eigenvalues of the matrix
$D_z g (y, h_0(y), 0)$
all have negative real parts.

\noindent\textbf{Remark.}
The two assumptions are justified
in most enzyme kinetics and
some combustion and atmospheric
chemistry problems.
In certain more complex
reaction mechanisms, however,
they may need justification.
Toward this end,
we observe that, 
in those cases where reduction
methods are expected to be effective,
$h_0$ can be found locally 
by the Implicit Function theorem
(since the second assumption
guarantees that the matrix
$(D_z g) (y, h_0(y), 0 )$
is invertible for each $y \in K$),
and GSPT can be applied to each
local portion.
In the absence of singularities,
these local functions
can be pieced together to form
a smooth global function
over the entire domain under consideration.~\qed

Under the above conditions,
standard asymptotic theory
(see, for example, Refs.~\cite{T-1948,
L-1950, F-1979, CH-1984, N-1985, MKKR-1994, J-1994})
guarantees that,
when $\eps$ is positive but
arbitrarily small,
there exists a \emph{slow manifold}
$\mathcal{M}_\eps$ 
that is invariant under the dynamics
of the system of Eqs.~(\ref{eq-y})--(\ref{eq-z}),
has the same dimension as $\mathcal{M}_0$,
and lies near $\mathcal{M}_0$.
All nearby solutions relax exponentially
fast to $\mathcal{M}_\eps$,
and their long-term evolution
is determined by an associated
solution on the slow manifold itself.
The manifold $\mathcal{M}_\eps$
is usually not unique;
there typically is a family
of slow manifolds, 
all exponentially close
($\mathcal{O}(\mathrm{e}^{-c/\eps})$ 
for some $c>0$).

\begin{THEOREM} \Label{t-Fenichel}
(Fenichel, asymptotically stable slow manifolds).
For any sufficiently small $\eps$,
there is a function $h_\eps$
that is defined on $K$
such that the graph
\be
  \mathcal{M}_\eps = \{ (y,z) : z = h_\eps (y) , \; y \in K \}
  \Label{M-eps}
\ee
is locally invariant 
under the dynamics of
Eqs.~(\ref{eq-y})--(\ref{eq-z}).
The function $h_\eps$ admits
an asymptotic expansion,
\be
  h_\eps (y)
  = h_0 (y) + \eps h^{(1)} (y)
  + \eps^2 h^{(2)} (y) + \cdots \quad
  \mbox{as } \eps \downarrow 0,
  \Label{h-eps-exp}
\ee
where the coefficients
$h^{(\ell)} : K \to \mathbf{R}^n$
are found successively from
the equation
\begin{eqnarray}
  (D_z g) h^{(\ell)}
  &\hspace{-0.5em}=\hspace{-0.5em}&
      \sum_{i=0}^{\ell-1} (Dh^{(i)}) f^{(\ell-1-i)} 
       - \sum_{j=2}^\ell \frac{1}{j!} (D_z^j g)
         \sum_{|i| = \ell} (h^{(i_1)}, \ldots \,, h^{(i_j)})  
      \nonumber \\
   &&\hspace{-2em}\mbox{}-
         \sum_{k=1}^{\ell-1} \frac{1}{k!} 
         \sum_{j=1}^{\ell-k} \frac{1}{j!} (D_z^j (\partial_\eps^k g))
                 \sum_{|i| = \ell-k} (h^{(i_1)}, \ldots\,, h^{(i_j)})
         - \frac{1}{\ell!} (\partial_\eps^\ell g)
  \Label{hl}
\end{eqnarray}
for $\ell = 1, 2, \ldots\,$,
with $h^{(0)} = h_0$.
Here, the functions 
$f$ and $g$ 
and their derivatives
are evaluated at $(y, z=h_0(y), 0)$,
and it is understood 
that a sum is empty 
when the lower bound 
exceeds the upper bound.
In particular,
$h^{(1)}$ and $h^{(2)}$
are given by
\begin{eqnarray}
  (D_z g) h^{(1)}
  &=& (Dh_0) f - g_\eps , \Label{h1} \\
  (D_z g) h^{(2)}
  &=& (D h^{(1)}) f
  + (Dh_0) \left( (D_z f) h^{(1)} + f_\eps \right) \nonumber \\
  &&\mbox{}-
  \half (D_z^2 g) \left( h^{(1)}, h^{(1)} \right)
  - (D_z g_\eps) h^{(1)} - \half g_{\eps\eps} .
\end{eqnarray}
Furthermore,
$h_\eps \in C^r (K)$
for any finite $r$,
and the dynamics 
of the system  of
Eqs.~(\ref{eq-y})--(\ref{eq-z})
on $\mathcal{M}_\eps$ 
are given by 
the reduced equation
\be
  \dot{y} = f ( y, h_\eps(y), \eps).
  \Label{traj}
\ee
\end{THEOREM} 

\begin{PROOF}
The theorem is 
a direct restatement 
of Theorem 2 in Ref.~\cite{J-1994}
for the special case
in which $\mathcal{M}_0$ 
is asymptotically stable.
It also follows directly
from~\cite[Theorem]{N-1985}
and is a special case
of the Fenichel theory~\cite{F-1979}.
The asymptotics of the
slow manifold $\mathcal{M}_\eps$
are given explicitly, for example,
in Refs.~\cite{N-1991, MKKR-1994}.
\end{PROOF}

\noindent\textbf{Remark.}
In many instances---for example,
in the Michaelis--Menten--Henri 
reaction mechanism
discussed in Section~\ref{s-MMH}
and various combustion problems---the
reduced slow system
${\dot y} = f (y, h_0(y), 0)$
has an asymptotically stable fixed point
at $(y_0, h_0 (y_0))$, say.
In such cases, the reaction scheme 
has a global attracting equilibrium.
Under the hypotheses made above, 
the system of
Eqs.~(\ref{eq-y})--(\ref{eq-z})
has a fixed point at
$(y_{0,\eps}, h_\eps (y_{0,\eps}))$,
and the slow manifold $\mathcal{M}_\eps$
is its weak stable manifold.~\qed

\noindent\textbf{Remark.}
While we have used it here only
for the case of attracting manifolds,
the Fenichel theorem 
and Theorem 2 in Ref. \cite{J-1994}
hold for the more general case 
of fast-slow systems of ODEs
for which the manifold $\mathcal{M}_0$
is normally hyperbolic---that is,
where there can be both
fast stable (exponentially contracting)
and fast unstable (exponentially expanding)
dynamics in the directions transverse 
to $\mathcal{M}_0$.
In the more general case,
the matrix $(D_z g) (y, h_0(y), 0)$
has $s$ eigenvalues with a negative real part
and $u$ eigenvalues with a positive real part,
the fast variable~$z$ decomposes into
a $u$-dimensional and a $s$-dimensional
component with $u+s=n$,
and the dynamics of all solutions
near $\mathcal{M}_\eps$
are governed by the
Fenichel normal
form~\cite{JKK-1996}.
The asymptotics of $\mathcal{M}_\eps$
remains unchanged.~\qed

\noindent\textbf{Remark.}
The articles of Tikhonov~\cite{T-1948}
and Levinson~\cite{L-1950, LL-1954}
present the original theory 
of persistence of asymptotically
stable manifolds
(see also Ref.~\cite{O-1991}).
The theory of persistence of
normally hyperbolic manifolds
can be found in the monographs
of Fenichel~\cite{F-1971,F-1979}
and Hirsch, Pugh, and Shub~\cite{HPS-1977};
see also Ref.~\cite{W-1994}.
Other relevant references are
\cite{N-1985} and~\cite{MKKR-1994}
for singularly perturbed systems
of ODEs with asymptotically stable slow manifolds
and~\cite{S-1990} and~\cite{J-1994}
for singularly perturbed systems of ODEs
with general normally hyperbolic slow manifolds.
An introductory exposition of GSPT
is given in Ref.~\cite{K-1999}.~\qed

\noindent\textbf{Remark.}
A numerical procedure for finding 
asymptotically stable slow manifolds
in fast-slow systems,
which is stable and highly accurate
for small values of $\eps$,
has been given by
Nipp~\cite{N-1991}.~\qed

\section{The ILDM Method of Maas and Pope \label{s-ILDM}}
\setcounter{equation}{0}
The ILDM method starts from the slow system,
Eqs.~(\ref{eq-y-slow})--(\ref{eq-z-slow}),
takes the local vector field~$F$
and the associated Jacobian~$J$,
and reduces the latter at each point
to a fast and a slow component.
The vector field $F$ and its
Jacobian $J$ are
\be
  F
  = \left(
  \begin{array}{c}
  f \\
  \eps^{-1} g
  \end{array}
  \right) , \quad
  J
  = \left(
  \begin{array}{cc}
  D_y f           & D_z f \\
  \eps^{-1} D_y g & \eps^{-1} D_z g
  \end{array}
  \right) ,
  \Label{F-jac}
\ee
where
$D_y f$ is the $m \times m$ matrix of
partial derivatives $\partial f_i / \partial y_j$,
$D_z f$ the $m \times n$ matrix of
partial derivatives $\partial f_i / \partial z_j$,
$D_y g$ the $n \times m$ matrix of
partial derivatives $\partial g_i / \partial y_j$,
and
$D_z g$ the $n \times n$ matrix of
partial derivatives $\partial g_i / \partial z_j$.

By assumption, the real part of each
eigenvalue of $J$ is negative.
The sum of the eigenvalues
is equal to the trace of $J$,
which is $\mathcal{O} (\eps^{-1})$
as $\eps \downarrow 0$,
and their product is equal to
the determinant of $J$,
which is $\mathcal{O} (\eps^{-n})$
as $\eps \downarrow 0$.
The eigenvalues of $J$ fall
therefore into two groups:
one group of $m$ eigenvalues
with $\mathcal{O} (1)$
negative real parts and
another group of $n$ eigenvalues
with $\mathcal{O} (\eps^{-1})$
negative real parts.
The eigenvectors associated with the
first group span the \emph{slow subspace},
those associated with the second group
the \emph{fast subspace}.
The Maas and Pope algorithm defines the
ILDM as the locus of all points $(y, z)$
where the vector field~$F$ lies
entirely in the slow subspace.

The algorithm uses a Schur
decomposition~\cite[Section~6.3]{S-1973}
of $J$,
\be
  J = Q N Q' ,
  \Label{schur-MP}
\ee
with $Q$ unitary
($QQ' = Q'Q = I_{m+n}$,
${}'$ denoting the transpose)
and $N$ upper triangular,
\be
  Q = (Q_s\ Q_f) , \quad
  N = \left( \begin{array}{cc}
  N_{s} & N_{sf} \\
  0     & N_f \end{array} \right) .
  \Label{QN-MP}
\ee
The dimensions of $Q_s$ and $Q_f$
are $(m+n) \times m$ and
$(m+n) \times n$, respectively;
$N_s$ is an $m \times m$ upper triangular
matrix,
$N_f$ an $n \times n$ upper triangular matrix,
and
$N_{sf}$ an $m \times n$ full matrix.
The eigenvalues of $J$
appear on the diagonal of $N$
in descending order
of their real parts,
from least negative
at the $(1,1)$ position
to most negative at the
$(m+n, m+n)$ position.
This particular ordering is accomplished
in Ref.~\cite{MP-1992}
by means of a modification of Stewart's
implementation of the Schur algorithm~\cite{S-1976}
and in Ref.~\cite{M-1998}
by means of a standard Schur decomposition
followed by 
a sequence of Givens
rotations~\cite[Section~5.1]{GL-1996}.

The first $m$ Schur vectors---that is,
the columns of $Q_s$---form
an orthogonal basis for the slow subspace,
while the remaining $n$ Schur
vectors---the columns of $Q_f$---form
an orthogonal basis for the
orthogonal complement of the slow subspace.
The vector field $F$ is entirely in
the slow subspace if it is orthogonal
to the orthogonal complement
of the slow subspace---that is,
if
\be
  Q'_f F = 0 .
  \Label{ILDM-MP}
\ee
This equation defines the ILDM,
the latter being an approximation of
the slow manifold~$\mathcal{M}_\eps$.
We analyze its asymptotics
(as $\eps \downarrow 0$)
in the following sections.

\noindent\textbf{Remark.}
The matrix $Q'_f$ corresponds to $Q_L^T$,
the number $n$ to $n_f$,
and the sum $m+n$ to $n$
in Ref.~\cite{MP-1992}.~\qed

In the numerical implementation
of the ILDM method for general,
closed, adiabatic, and isobaric
reaction mechanisms,
the system of equations is closed
by supplementing the ILDM equation,
Eq.~(\ref{ILDM-MP}),
by a set of parameter equations.
The parameter equations fix the enthalpy,
pressure, and element composition.
In addition, the reaction progress variables
are treated as parameters.
Each fixed set of parameters yields
one point of the ILDM,
and the entire ILDM is obtained
by sweeping over the admissible
set of parameter values.
As noted in Ref.~\cite{MP-1992},
the parameters can generally be chosen
so the ILDM is at least defined piecewise,
and, most important, the choice of
the parameter equations does not influence
the construction of the manifold.

In Eqs.~(\ref{eq-y})--(\ref{eq-z}),
the enthalpy, pressure, and conserved
quantities have been neglected.
In this case, the parameter equations
fix the values of the slow variables $y$,
and the ILDM is obtained by sweeping
over all points $y \in K$.

\section{Asymptotics of the ILDM --- Planar Case \label{s-ILDM-planar}}
\setcounter{equation}{0}
We first restrict our attention
to planar fast-slow systems,
Eqs.~(\ref{eq-y-slow})--(\ref{eq-z-slow})
with $m = n = 1$,
for which the computations
are relatively straightforward
and the asymptotic analysis
more transparent.
We address the general case
in Section~\ref{s-ILDM-general}.

In the planar case,
the vector field $F$
and its Jacobian $J$
are
\be
  F
  = \left( 
  \begin{array}{c}
  f \\ \eps^{-1} g \end{array}
  \right) , \quad
  J
  = \left(
  \begin{array}{cc}
  f_y           & f_z \\
  \eps^{-1} g_y & \eps^{-1} g_z
  \end{array}
  \right) .
  \Label{jac-planar}
\ee
The eigenvalues of $J$ are
\be
  \lambda_{s,f}
  =\mbox{} \half
  \left( \eps^{-1} g_z + f_y \right)
  \pm \sqrt{ \onefourth \left(\eps^{-1} g_z + f_y \right)^2
  - \eps^{-1} \left( f_y g_z - f_z g_y \right)} ,
  \Label{lambdapm}
\ee
where the upper (lower) sign
is associated with
$\lambda_s$ ($\lambda_f$).
Thus,
\be
  \lambda_s = f_y - \frac{f_z g_y}{g_z} + \mathcal{O} (\eps) , \quad
  \lambda_f = \eps^{-1} g_z + \mathcal{O} (1) \quad \mbox{as } \eps \downarrow 0 .
  \Label{lambda-sf-as}
\ee
The derivatives of $f$ and $g$,
which are evaluated at $(y, z, \eps)$,
are all $\mathcal{O} (1)$ as $\eps \downarrow 0$.
The (nonnormalized) slow eigenvector is
\be
  v_{s} = \left(
  \begin{array}{c}
  \lambda_{s} - \eps^{-1} g_z \\ \eps^{-1} g_y
  \end{array}
  \right) ,
  \Label{vs}
\ee
and there is a corresponding fast eigenvector~$v_f$.
The vector $v_s$ spans the slow subspace,
$v_f$ the fast subspace.
The vectors $v_s$ and $v_f$ are
not necessarily orthogonal.
To determine the points $(y, z)$
in the phase plane where the vector field $F$
lies entirely in the slow subspace,
we work with the orthogonal complement
of the slow subspace,
which is spanned by
the row vector
\be
  v_s^\perp
  = ( \eps^{-1} g_y , \, \eps^{-1} g_z - \lambda_{s}) .
  \Label{vs-perp}
\ee
The locus of all points
in the phase plane where
the vector field $F$
is in the slow subspace
coincides with the set of
all points $(y, z)$
where $F$ is orthogonal to
$v_s^\perp$---that is, where
\be
  f g_y + g (\eps^{-1} g_z - \lambda_{s}) = 0 .
  \Label{ILDM-planar}
\ee
This equation defines the ILDM.

\noindent\textbf{Remark.}
A Schur decomposition
of the matrix $J$
gives the vectors
$q_{f} = v_f / |v_f|$
and
$q_{f}^\perp = v_f^\perp / |v_f|$
directly.
The algorithm must be modified
to find the vectors
$q_{s} = v_s / |v_s|$
and
$q_{s}^\perp = v_s^\perp / |v_s|$,
as described in Section~\ref{s-ILDM}.~\qed

\begin{THEOREM} \Label{t-ILDM-planar}
(Planar Case).
The equation for the ILDM,
Eq.~(\ref{ILDM-planar}), admits
an asymptotic solution in the form
of a power series expansion,
\be
  z = \psi (y, \eps)
  = \psi^{(0)} (y)
  + \eps \psi^{(1)} (y)
  + \eps^2 \psi^{(2)} (y)
  + \cdots
  \quad \mbox{as } \eps \downarrow 0 .
  \Label{z-exp-planar}
\ee
The functions
$\psi^{(0)}$, $\psi^{(1)}$, and $\psi^{(2)}$
are defined by the equations
\begin{eqnarray}
  \psi^{(0)}
  &\hspace{-0.5em}=\hspace{-0.5em}&
  h_0 , \Label{MP-z0} \\
  g_z \psi^{(1)}
  &\hspace{-0.5em}=\hspace{-0.5em}&
  f h_0' - g_\eps , \Label{MP-z1} \\
  g_z \psi^{(2)}
  &\hspace{-0.5em}=\hspace{-0.5em}& 
     f \psi^{(1)'}
       - \frac{f^2}{g_z} h_0''
       + (f_z \psi^{(1)} + f_\eps) h_0' \nonumber \\
  && \mbox{}-
       \half g_{zz} (\psi^{(1)})^2
       - g_{z\eps} \psi^{(1)} 
       - \half g_{\eps\eps} .  \Label{MP-z2}
\end{eqnarray}
Here, $h_0 \equiv h_0(y)$
is defined by the equation
$g ( y, h_0 (y), 0) = 0$,
${}'$ denotes differentiation
with respect to $y$,
and the functions $f$ and $g$
and their derivatives 
are evaluated 
at $(y, h_0 (y), 0)$.
\end{THEOREM}

\begin{PROOF}
Assume that $z = \psi (y, \eps)$,
where $\psi$ is given by the power
series expansion~(\ref{z-exp-planar}).
Then
\begin{eqnarray}
\lefteqn{ f (y, \psi(y, \eps), \eps)
  = f
  + \eps \left( f_{z} \psi^{(1)}
  + f_{\eps} \right) } &&\mbox{} \nonumber \\
  &&\mbox{} \hspace{4em}
  + \eps^2 \left(f_{z} \psi^{(2)} + \half f_{zz} (\psi^{(1)})^2
            + f_{z\eps} \psi^{(1)} + \half f_{\eps\eps} \right)
  + \cdots \,,
  \Label{f-exp-planar}
\end{eqnarray}
where, in the right member,
$f$ and its derivatives
are evaluated at
$(y, \psi^{(0)} (y), 0)$.
Similar expansions hold for $g$
and the derivatives of $f$ and $g$.
The leading term in the expansion
of $\lambda_{s}$ follows immediately
from Eq.~(\ref{lambda-sf-as}),
\be
  \lambda_{s}
  = \lambda^{(0)}_{s}
    + \mathcal{O}(\eps) , \quad
  \lambda^{(0)}_{s}
  = f_y - \frac{f_z g_y}{g_z} ,
   \Label{lambda0s} 
\ee
where the derivatives of $f$ and $g$ are
similarly evaluated at $(y, \psi^{(0)} (y), 0)$.
We substitute the various expansions 
into Eq.~(\ref{ILDM-planar})
and equate the coefficients
of like powers of $\eps$.

\noindent
$\mathbf{\mathcal{O}(\eps^{-1})}$.
The ILDM equation, Eq.~(\ref{ILDM-planar}),
gives
\be
  (g g_z) (y, z = \psi^{(0)} (y), 0) = 0 ,
  \Label{eq0}
\ee
which is satisfied if
$\psi^{(0)} = h_0$.
This result confirms Eq.~(\ref{MP-z0}).

\noindent
$\mathbf{\mathcal{O}(1)}$.
From the equation for the ILDM,
Eq.~(\ref{ILDM-planar}), we obtain
\be
  f g_y + (g_z \psi^{(1)} + g_\eps) g_z = 0 .
  \Label{eq1}
\ee
Here,
we have used 
the identity
$g \equiv g (y, h_0 (y), 0) = 0$.
The same identity implies that
\be
  g_y + g_z h_0' = 0 ,
  \Label{gy-planar}
\ee
so Eq.~(\ref{eq1}) reduces to
\be
  \left( g_z \psi^{(1)} + g_\eps - f h_0' \right) g_z = 0 .
  \Label{gz}
\ee
The assumption of attractive manifolds
implies that $g_z < 0$,
so Eq.~(\ref{MP-z1}) follows.

\noindent
$\mathbf{\mathcal{O}(\eps)}$.
From the equation for the ILDM,
Eq.~(\ref{ILDM-planar}), we obtain
\begin{eqnarray}
  &&
  f \left( g_{yz} \psi^{(1)} + g_{y\eps}
  + (g_{zz} \psi^{(1)} + g_{z\eps} - \lambda^{(0)}_s) h_0' \right)
  - (f_z \psi^{(1)} + f_\eps) g_z h_0' \nonumber \\
  &&\mbox{}
  + (g_z \psi^{(2)} + \half g_{zz} (\psi^{(1)})^2
        + g_{z\eps} \psi^{(1)} + \half g_{\eps\eps}) g_z
  = 0 .
  \Label{eq2}
\end{eqnarray}
Here, we have used
Eqs.~(\ref{MP-z0})
and~(\ref{MP-z1})
and the identity~(\ref{gy-planar}).
The same identity
also results in a
simplification of
the expression~(\ref{lambda0s})
for $\lambda^{(0)}_{s}$,
\be
 \lambda^{(0)}_{s}
 = f_y + f_z h_0' .
 \Label{lambda0s-red}
\ee
Furthermore,
differentiating Eq.~(\ref{gz})
with respect to $y$,
we find
\be
  g_{yz} \psi^{(1)} + g_{y\eps }
  + (g_{zz} \psi^{(1)} + g_{z\eps}) h_0'
  + g_z \psi^{(1)'}
  = f h_0'' + (f_y + f_z h_0') h_0' .
  \Label{gyz-planar}
\ee
With Eqs.~(\ref{lambda0s-red})
and~(\ref{gyz-planar}),
Eq.~(\ref{eq2}) simplifies to
\begin{eqnarray}
\lefteqn{
  \left( \rule{0em}{2ex} g_z \psi^{(2)}
  + \half g_{zz} (\psi^{(1)})^2
  + g_{z\eps} \psi^{(1)}
  + \half g_{\eps\eps}
  \right. }&&\mbox{} \nonumber \\
  &&\left.\mbox{}-
  f \psi^{(1)'}
  + (f^2/g_z) h_0''
  - (f_z \psi^{(1)} + f_\eps) h_0' \right) g_z = 0 .
\end{eqnarray}
Since $g_z < 0$, Eq.~(\ref{MP-z2}) follows.
\end{PROOF}

In the following section, we will
generalize Theorem~\ref{t-ILDM-planar}
to the multidimensional case
(Theorem~\ref{t-ILDM-general})
and compare the asymptotics of
the ILDM with the asymptotics of
the slow manifold $\mathcal{M}_\eps$
(Corollary~\ref{c-ILDM}).

\section{Asymptotics of the ILDM --- General Case \label{s-ILDM-general}}
\setcounter{equation}{0}
The definition of the ILDM,
Eq.~(\ref{ILDM-MP}),
is based on a partition
of the Jacobian,
Eq.~(\ref{schur-MP}),
into a fast and a slow component
at each point of phase space
and a Schur decomposition
to generate bases for the
corresponding fast and slow subspaces.
Practical implementations of the
Schur decomposition rely typically
on the method of deflation~\cite[Chapter~7]{GL-1996};
hence, the eigenvalues are generated
in the order of descending
\emph{absolute values}
of their real parts.
This procedure yields
a unitary matrix of the form
$Q = (Q_f\ Q_s)$.
The columns of $Q$ are then reordered,
for example by a sequence of
Givens rotations,
as in Ref.~\cite{M-1998}.

Although this procedure is practical
for numerical computations,
it is not amenable to analysis.
We start therefore
from the standard
Schur decomposition
\emph{before} reordering,
\be
  J = QTQ' ,
  \Label{schur}
\ee
where
\be
  T = \left( \begin{array}{cc}
  \Lambda_{f} & \Lambda \\
  0           & \Lambda_{s} \end{array} \right) ,
  \Label{T-decomp}
\ee
with
$\Lambda_{f}$ an $n \times n$ upper triangular matrix,
$\Lambda_{s}$ an $m \times m$ upper triangular matrix,
and $\Lambda$ an $n \times m$ full matrix.
The diagonal elements of $\Lambda_{f}$
are the $\mathcal{O}(\eps^{-1})$ eigenvalues of $J$,
and the diagonal elements of $\Lambda_{s}$
are the $\mathcal{O}(1)$ eigenvalues of $J$.
The structure of the unitary matrix $Q$ is
\be
  Q = \left( \begin{array}{cc}
            Q_{11} & Q_{12} \\
            Q_{21} & Q_{22} \end{array} \right) ,
  \Label{Q-decomp}
\ee
where
$Q_{11}$ is an $m \times n$ matrix,
$Q_{12}$ an $m \times m$ matrix,
$Q_{21}$ an $n \times n$ matrix, and
$Q_{22}$ an $n \times m$ matrix.
The columns of
$\left( \begin{array}{c} Q_{11} \\ Q_{21} \end{array} \right)$
and
$\left( \begin{array}{c} Q_{12} \\ Q_{22} \end{array} \right)$
form an orthogonal basis of the fast subspace
and its orthogonal complement, respectively.

Since the fast and slow subspaces are not necessarily
mutually orthogonal, the orthogonal complement
of the fast subspace does not necessarily coincide
with the slow subspace, and a further operation
is needed to identify a basis for the slow subspace.
This operation consists of solving
the Sylvester equation
\be
  \Lambda_{f} X - X \Lambda_{s} = - \Lambda
  \Label{sylv}
\ee
for the $n \times m$ matrix $X$.
With the definition
\be
  Y = \left( \begin{array}{cc} I_n & X \\ 0 & I_m \end{array} \right) ,
\ee
we obtain a \emph{block diagonalization} of $J$,
\be
  J = (QY) T_d (QY)^{-1} ,
  \Label{decomp}
\ee
where
\be
  T_d = \left( \begin{array}{cc}
  \Lambda_{f} & 0 \\
  0           & \Lambda_{s} \end{array} \right) ,
  \Label{Td}
\ee
\be
  QY = \left( \begin{array}{cc} Q_{11} & Q_{11} X + Q_{12} \\
                                Q_{21} & Q_{21} X + Q_{22} \end{array} \right) , \quad
  (QY)^{-1} = \left( \begin{array}{cc} Q'_{11} - X Q'_{12} & Q'_{21} - X Q'_{22} \\
                                Q'_{12} & Q'_{22} \end{array} \right) .
  \Label{QY}
\ee
Thus, $QY$ reduces the matrix $J$ to its fast
and slow components, and the condition that
the vector field $F$ given in Eq.~(\ref{F-jac})
must lie entirely in the slow subspace
is satisfied if
\be
  (Q'_{11} - X Q'_{12}) f + \eps^{-1} (Q'_{21} - X Q'_{22}) g = 0 .
  \Label{ILDM}
\ee
The ILDM obtained from Eq.~(\ref{ILDM})
is the same as the ILDM obtained from
Eq.~(\ref{ILDM-MP}) and also
the same as the ILDMs obtained in
Refs.~\cite{MP-1992} and~\cite{M-1998}.

\begin{THEOREM} \Label{t-ILDM-general}
(General Case).
The equation for the ILDM,
Eq.~(\ref{ILDM-MP}), admits
an asymptotic solution in the form
of a power series expansion,
\be
  z = \psi (y, \eps)
  = \psi^{(0)} (y)
  + \eps \psi^{(1)} (y)
  + \eps^2 \psi^{(2)} (y)
  + \cdots
  \quad \mbox{as } \eps \downarrow 0 .
  \Label{psi-exp}
\ee
The $\mathbf{R}^n$-valued functions
$\psi^{(0)}$, $\psi^{(1)}$, and $\psi^{(2)}$
are defined by the equations
\begin{eqnarray}
  \psi^{(0)} &=& h_0 ,  \Label{MP-0} \\
  (D_z g) \psi^{(1)} &=&
  (Dh_0) f - g_\eps , \Label{MP-1} \\
  (D_z g) \psi^{(2)} &=&
  (D \psi^{(1)}) f - (D_z g)^{-1} (D^2h_0) (f, f)
  + (Dh_0) \left((D_z f) \psi^{(1)} + f_\eps \right) \nonumber \\
  &&\mbox{}- \half (D^2_z g) (\psi^{(1)}, \psi^{(1)})
  - (D_z g_\eps) \psi^{(1)}
  - \half g_{\eps\eps} .
  \Label{MP-2}
\end{eqnarray}
Here, $h_0 \equiv h_0 (y)$ is
the $\mathbf{R}^n$-valued function
defined by Eq.~(\ref{h0}),
$g (y, h_0(y), 0) = 0$;
$Dh_0 \equiv (Dh_0) (y)$ is
a linear operator from
$\mathbf{R}^m$ to $\mathbf{R}^n$,
which is represented by the
$n \times m$ matrix of partial derivatives
$\partial h_{0,i} / \partial y_j$,
and $D^2h_0 = D(Dh_0) \equiv (D^2h_0) (y)$ is
a bilinear map from
$\mathbf{R}^m \times \mathbf{R}^m$
to $\mathbf{R}^n$,
$(D^2h_0) (u, v) = ((D^2h_0)u)v$
for all $u, v \in \mathbf{R}^m$.
The functions $f$ and $g$
and their derivatives 
are evaluated 
at $(y, h_0 (y), 0)$;
$D_z f \equiv D_z f (y, h_0 (y), 0)$
is a linear operator from
$\mathbf{R}^n$ to $\mathbf{R}^m$,
$D_z g \equiv D_z g (y, h_0 (y), 0)$
a linear operator from
$\mathbf{R}^n$ to $\mathbf{R}^n$, and
$D_z^2 g = D_z (D_z g) \equiv D_z^2 g (y, h_{0} (y), 0)$
a bilinear map from
$\mathbf{R}^n \times \mathbf{R}^n$ to $\mathbf{R}^n$.
\end{THEOREM}

\begin{PROOF}
Assume that $z = \psi (y, \eps)$
and that $\psi$ is given by the
expansion~(\ref{psi-exp}).
For the asymptotic analysis
of Eq.~(\ref{ILDM}),
we take
\be
  Q \equiv Q (\eps)
  =
  \left( \begin{array}{cc}
            0            & I_m \\
            Q^{(0)}_{21} & 0   \end{array}
  \right)
  + \eps
  \left( \begin{array}{cc}
            Q^{(1)}_{11} & 0            \\
            0            & - Q^{(0)}_{21} Q^{(1)'}_{11} \end{array}
  \right)
  + \cdots \,,
  \Label{Q-exp}
\ee
with $Q^{(0)}_{21}$ a unitary $n \times n$ matrix
and $Q^{(1)}_{11}$ an $m \times n$ matrix
to be determined.
Thus, $Q$ is unitary to $\mathcal{O} (\eps)$.
Higher-order terms can be found in a consistent manner
so $Q(\eps)$ is unitary to any desired order.
We take, furthermore,
\begin{eqnarray}
  \Lambda_{f} &\equiv& \Lambda_{f} (\eps)
  = \eps^{-1} \Lambda_{f}^{(-1)}
  + \Lambda_{f}^{(0)} + \cdots \,,
  \Label{Lm-exp} \\
  \Lambda &\equiv& \Lambda (\eps)
  = \eps^{-1} \Lambda^{(-1)}
  + \Lambda^{(0)} + \cdots \,,
  \Label{L-exp} \\
  \Lambda_{s} &\equiv& \Lambda_{s} (\eps)
  = \Lambda_{s}^{(0)} + \cdots \,,
  \Label{Lp-exp} \\
  X &\equiv& X (\eps)
  = X^{(0)} + \eps X^{(1)} + \cdots \,.
  \Label{X-exp}
\end{eqnarray}
The generalization of
the expansion~(\ref{f-exp-planar})
to the present case is
\begin{eqnarray}
\lefteqn{ f (y, \psi(y, \eps), \eps)
  = f
  + \eps \left( (D_z f) \psi^{(1)} + f_{\eps} \right)
  + \eps^2 \left( (D_z f) \psi^{(2)} \right. }&&\mbox \nonumber \\
  &&\hspace{5em} \left. \mbox{}
  + \half (D^2_z f) (\psi^{(1)}, \psi^{(1)})
            + (D_z f_{\eps}) \psi^{(1)} + \half f_{\eps\eps} \right)
  + \cdots \,.
  \Label{f-exp-gen}
\end{eqnarray}
In the right member,
$f$ and its derivatives are
evaluated at $(y, \psi^{(0)} (y), 0)$.
Similar expansions hold for $g$
and the derivatives of $f$ and $g$.

To prove the theorem,
we substitute the various expansions
into Eq.~(\ref{ILDM})
and equate the coefficients
of like powers in $\eps$
in the usual manner.

\noindent
$\mathbf{\mathcal{O}(\eps^{-1})}$.
The ILDM equation, Eq.~(\ref{ILDM}), gives
\be
  Q^{(0)'}_{21} g = 0 .
  \Label{ILDM-0}
\ee
Since $Q^{(0)}_{21}$ is unitary,
Eq.~(\ref{ILDM-0}) reduces to
\be
  g (y, \psi^{(0)} (y), 0) = 0 .
  \Label{g}
\ee
This equation is satisfied if $\psi^{(0)} = h_0$,
which confirms Eq.~(\ref{MP-0}).

\noindent
$\mathbf{\mathcal{O}(1)}$.
From the ILDM equation, Eq.~(\ref{ILDM}),
we obtain
\be
  - X^{(0)} f
  + Q^{(0)'}_{21} \left( (D_z g) \psi^{(1)} + g_\eps \right)
  = 0 .
  \Label{ILDM-1}
\ee
Here, we have already used the
identity~$g = 0$.

The matrix $X^{(0)}$ is determined from
the $\mathcal{O} (\eps^{-1})$ terms in
the Sylvester equation, Eq.~(\ref{sylv}),
\be
  \Lambda^{(-1)}_{f} X^{(0)} = - \Lambda^{(-1)} .
  \Label{sylv-0}
\ee
The matrices $\Lambda^{(-1)}_{f}$ and $\Lambda^{(-1)}$,
in turn, follow from
the $\mathcal{O} (\eps^{-1})$ terms in
the Schur decomposition, Eq.~(\ref{schur}),
\be
  Q^{(0)}_{21} \Lambda^{(-1)}
  = D_y g , \quad
  Q^{(0)}_{21} \Lambda^{(-1)}_{f} Q^{(0)'}_{21}
  = D_z g .
  \Label{eq-067-red}
\ee
The second equation is the
Schur decomposition of $D_z g$,
so $Q^{(0)}_{21}$ is determined
by the ordering of the elements
of $\Lambda^{(-1)}_{f}$.
Both equations can be inverted,
\be
  \Lambda^{(-1)} = Q^{(0)'}_{21} (D_y g) , \quad
  \Lambda^{(-1)}_{f} = Q^{(0)'}_{21} (D_z g) Q^{(0)}_{21} .
  \Label{lambda-1}
\ee
Hence,
\be
  X^{(0)}
  = - Q^{(0)'}_{21} (D_z g)^{-1} (D_y g) .
  \Label{X0}
\ee
We can simplify this expression
if we use the identity
$g (y, h_0(y), 0) = 0$,
which holds for all $y$.
Upon differentiation,
the identity gives
a relation between
$D_y g$ and $D_z g$,
\be
  D_y g + (D_z g) (Dh_0) = 0 .
  \Label{gy}
\ee
Note that this is a relation
in the space of linear operators
from $\mathbf{R}^m$ to $\mathbf{R}^n$.
With this identity,
Eq.~(\ref{X0}) becomes
\be
  X^{(0)} = Q^{(0)'}_{21} (Dh_0) ,
  \Label{X0-red}
\ee 
and Eq.~(\ref{ILDM-1}) reduces to
\be
  Q^{(0)'}_{21}
  \left[ (D_z g) \psi^{(1)} + g_\eps - (Dh_0) f \right] = 0 .
\ee
Since $Q^{(0)}_{21}$ is unitary,
Eq.~(\ref{MP-1}) follows.

\noindent
$\mathbf{\mathcal{O}(\eps)}$.
From the ILDM equation, Eq.~(\ref{ILDM}),
we obtain
\begin{eqnarray}
  &&\left(
  Q^{(1)'}_{11} - X^{(1)}
  + Q^{(0)'}_{21} (Dh_0) Q^{(1)}_{11} Q^{(0)'}_{21} (Dh_0)
  \right) f
  - Q^{(0)'}_{21} (Dh_0) \left((D_z f) \psi^{(1)} + f_\eps\right) \nonumber \\
  &&\mbox{}+ Q^{(0)'}_{21}
  \left( (D_z g) \psi^{(2)}
  + \half (D_z^2 g) (\psi^{(1)}, \psi^{(1)}) + (D_z g_\eps) \psi^{(1)}
  + \half g_{\eps\eps} \right)
  = 0 . \Label{ILDM-2}
\end{eqnarray}
Here, we have already made use
of Eqs.~(\ref{MP-0}) and~(\ref{MP-1})
and substituted
the expression~(\ref{X0-red})
for $X^{(0)}$.

The matrix $X^{(1)}$ is determined from
the $\mathcal{O} (1)$ terms in
the Sylvester equation, Eq.~(\ref{sylv}),
\be
  \Lambda^{(-1)}_{f} X^{(1)} + \Lambda^{(0)}_{f} X^{(0)}
  - X^{(0)} \Lambda^{(0)}_{s} = - \Lambda^{(0)} .
  \Label{sylv-1}
\ee
The matrices $\Lambda^{(0)}_{f}$,
$\Lambda^{(0)}_{s}$, and $\Lambda^{(0)}$
follow in turn from the
$\mathcal{O} (1)$ terms in
the Schur decomposition, Eq.~(\ref{schur}),
\begin{eqnarray}
  Q^{(1)}_{11} \Lambda^{(-1)}
  + \Lambda^{(0)}_{s} &=& D_y f , \Label{eq-14} \\
  Q^{(1)}_{11} \Lambda^{(-1)}_{f} Q^{(0)'}_{21} &=& D_z f  , \Label{eq-15} \\
  Q^{(0)}_{21} \left(\Lambda^{(-1)}_{f} Q^{(1)'}_{11}
  + \Lambda^{(0)} \right)
  &=&
  (D_z (D_y g)) \psi^{(1)} + D_y g_{\eps} , \Label{eq-16} \\
  Q^{(0)}_{21} \left( - \Lambda^{(-1)} Q^{(1)}_{11} Q^{(0)'}_{21}
  + \Lambda^{(0)}_{f} Q^{(0)'}_{21} \right)
  &=&
  (D_z^2 g) \psi^{(1)} + D_z g_{\eps} . \Label{eq-17}
\end{eqnarray}
We proceed as follows.
First, we solve Eq.~(\ref{eq-15})
for $Q^{(1)}_{11}$,
\be
  Q^{(1)}_{11}
  = (D_z f) Q^{(0)}_{21} \left(\Lambda^{(-1)}_{f}\right)^{-1}
  = (D_z f) (D_z g)^{-1} Q^{(0)}_{21} .
  \Label{q111}
\ee
Then, we obtain $\Lambda^{(0)}_{s}$ from Eq.~(\ref{eq-14}),
\be
  \Lambda^{(0)}_{s}
  =
  D_y f - Q^{(1)}_{11} \Lambda^{(-1)}
  =
  D_y f + (D_z f) (Dh_0) .
\ee
(We have used the relation~(\ref{gy})
to rewrite the expression~(\ref{lambda-1})
for $\Lambda^{(-1)}$.)
Next, we solve
Eqs.~(\ref{eq-16}) and~(\ref{eq-17})
for $\Lambda^{(0)}$ and $\Lambda^{(0)}_{f}$,
\begin{eqnarray}
  \Lambda^{(0)}
  &\hspace{-0.5em}=\hspace{-0.5em}&
  Q^{(0)'}_{21}
  \left( (D_z (D_y g)) \psi^{(1)} + D_y g_{\eps}
  - Q^{(0)}_{21} \Lambda^{(-1)}_{f} Q^{(1)'}_{11} \right) \nonumber \\
  &\hspace{-0.5em}=\hspace{-0.5em}&
  Q^{(0)'}_{21}
  \left(
  (D_z (D_y g)) \psi^{(1)} + D_y g_{\eps}
  - (D_z g) \left((D_z f) (D_z g)^{-1}\right)'
  \right) , \\
  \Lambda^{(0)}_{f}
  &\hspace{-0.5em}=\hspace{-0.5em}&
  Q^{(0)'}_{21}
  \left( (D_z^2 g) \psi^{(1)} + D_z g_{\eps}
  + Q^{(0)}_{21} \Lambda^{(-1)}
    Q^{(1)}_{11} Q^{(0)'}_{21} \right) Q^{(0)}_{21} \nonumber \\ 
  &\hspace{-0.5em}=\hspace{-0.5em}&
  Q^{(0)'}_{21}
  \left(
  (D_z^2 g) \psi^{(1)} + D_z g_{\eps}
  - (D_z g) (Dh_0) (D_z f) (D_z g)^{-1}
  \right)
  Q^{(0)}_{21} .
\end{eqnarray}
After these steps,
we find $X^{(1)}$
from Eq.~(\ref{sylv-1}),
\begin{eqnarray}
  X^{(1)}
  &=&
  ( \Lambda^{(-1)}_{f} )^{-1}
  \left( - \Lambda^{(0)}
         - \Lambda^{(0)}_{f} X^{(0)}
         + X^{(0)} \Lambda^{(0)}_{s}
  \right) \nonumber \\
  &=&
  Q^{(0)'}_{21} (D_z g)^{-1}
  \left[
  - \left( (D_z (D_y g)) \psi^{(1)} + D_y g_\eps \right)
  + (D_z g) \left( (D_z f) (D_z g)^{-1} \right)'
  \right. \nonumber \\
  &&\mbox{}
  - \left( (D_z^2 g) \psi^{(1)} + D_z g_\eps
        - (D_z g)(Dh_0) (D_z f) (D_z g)^{-1} \right) (Dh_0) \nonumber \\
  &&\left.\mbox{}
  + (Dh_0) \left( D_y f + (D_z f)(Dh_0) \right)
  \rule{0em}{3ex}\right] .
  \Label{x1}
\end{eqnarray}
Substituting
$Q_{11}^{(1)}$ from Eq.~(\ref{q111})
and $X^{(1)}$ from Eq.~(\ref{x1})
into Eq.~(\ref{ILDM-2}),
we obtain
\begin{eqnarray}
  &Q^{(0)'}_{21} \hspace{-0.5em}&
  \left[ \rule{0em}{2ex}
  (D_z g)^{-1}
  \left(
  (D_z (D_y g)) (\psi^{(1)}, f)
  + (D_y g_\eps) f
  + (D_z^2 g) (\psi^{(1)}, (Dh_0) f)
  \right.
  \right. \nonumber \\
  &&\left.\hspace{3em}\mbox{}+
  (D_z g_\eps) ((Dh_0) f)
  - (Dh_0) \left( D_y f + (D_z f)(Dh_0) \right) f
  \rule{0em}{2ex}\right) \nonumber \\
  &&\mbox{}-
  (Dh_0) \left((D_z f) \psi^{(1)} + f_\eps\right) \nonumber \\
  &&\left.\mbox{}+
  (D_z g) \psi^{(2)}
  + \half (D_z^2 g) (\psi^{(1)}, \psi^{(1)})
  + (D_z g_\eps) \psi^{(1)}
  + \half g_{\eps\eps}
  \rule{0em}{2ex}\right]
  = 0 .
  \Label{ILDM-2-full}
\end{eqnarray}
The bilinear maps
$D_z (D_y g)$ and $D_z^2 g$
satisfy the symmetry relations
\begin{eqnarray}
  (D_z (D_y g)) (u, v) &=& (D_y (D_z g)) (v, u) , \quad
  u \in \mathbf{R}^n, \, v \in \mathbf{R}^m , \\
  (D_z^2 g) (u, v) &=& (D_z^2 g) (v, u) , \quad
  u, v \in \mathbf{R}^n ,
\end{eqnarray}
so Eq.~(\ref{ILDM-2-full})
is equivalent with
\begin{eqnarray}
  &Q^{(0)'}_{21} \hspace{-0.5em}&
  \left[ \rule{0em}{2ex}
  (D_z g)^{-1}
  \left(
  (D_y (D_z g)) (f, \psi^{(1)})
  + (D_y g_\eps) f
  + (D_z^2 g) ((Dh_0) f, \psi^{(1)})
  \right. \right. \nonumber \\
  &&\left.\hspace{3em}\mbox{}+
  (D_z g_\eps) ((Dh_0) f)
  - (Dh_0) \left( D_y f + (D_z f)(Dh_0) \right) f
  \rule{0em}{2ex}\right) \nonumber \\
  &&\mbox{}-
  (Dh_0) \left((D_z f) \psi^{(1)} + f_\eps\right) \nonumber \\
  &&\left.\mbox{}+
  (D_z g) \psi^{(2)}
  + \half (D_z^2 g) (\psi^{(1)}, \psi^{(1)})
  + (D_z g_\eps) \psi^{(1)}
  + \half g_{\eps\eps}
  \rule{0em}{2ex}\right]
  = 0 . \Label{ILDM-2-rev}
\end{eqnarray}
We simplify this expression
by means of Eq.~(\ref{MP-1}).
Upon differentiation,
this equation gives the identity
\begin{eqnarray}
  (D_y (D_z g)) \psi^{(1)}
  &\hspace{-0.5em}+\hspace{-0.5em}& D_y g_\eps
  + (D_z^2 g) (Dh_0) \psi^{(1)} + (D_z g_\eps)(Dh_0)
  + (D_z g) (D\psi^{(1)}) \nonumber \\
  &&\mbox{}= (D^2 h_0) f
  + (Dh_0) \left( D_y f + (D_z f) (Dh_0) \right) .
  \Label{gyz}
\end{eqnarray}
This is a relation in the space of
linear operators from $\mathbf{R}^m$
to $\mathbf{R}^n$.
When applied to the vector $f$,
it gives the identity
\begin{eqnarray}
  &(D_y (D_z g)) (f, \psi^{(1)})
  + (D_y g_\eps) f
  + (D_z^2 g) ((Dh_0)f, \psi^{(1)})
  + (D_z g_\eps) ((Dh_0)f) 
  \hspace{2em}&\mbox{} \nonumber \\
  &\mbox{}+
  (D_z g) ((D\psi^{(1)}) f)
  = (D^2 h_0) (f, f)
  + (Dh_0) \left( D_y f + (D_z f) (Dh_0) \right) f .
  \hspace{2em}&\mbox{}
  \Label{gyzf}
\end{eqnarray}
With this result,
Eq.~(\ref{ILDM-2-rev}) simplifies to
\begin{eqnarray}
  &Q^{(0)'}_{21}
  \left[
  (D_z g) \psi^{(2)}
  + \half (D^2_z g) (\psi^{(1)}, \psi^{(1)})
  + (D_z g_\eps) \psi^{(1)}
  + \half g_{\eps\eps}
  \right. 
  \hspace{2em}& \mbox{} \nonumber \\
  &\left.\mbox{}-
  (D\psi^{(1)}) f
  + (D_z g)^{-1} (D^2 h_0) (f, f)
  - (Dh_0) \left( (D_z f) \psi^{(1)} + f_\eps \right)
  \right]
  = 0 . 
  \hspace{2em}& \mbox{}
\end{eqnarray}
Since $Q^{(0)}_{21}$ is unitary,
Eq.~(\ref{MP-2}) follows.
\end{PROOF}

The following corollary summarizes
the result of the asymptotic analysis.

\begin{COROLLARY}  \Label{c-ILDM}
The ILDM is an approximation to the
slow manifold $\mathcal{M}_\eps$
of the fast-slow system of
Eqs.~(\ref{eq-y})--(\ref{eq-z}),
which is asymptotically accurate
up to and including the
$\mathcal{O} (\eps)$ term
as $\eps \downarrow 0$.
The approximation is
asymptotically accurate
up to and including the
$\mathcal{O} (\eps^2)$ term
if and only if
$D^2 h_0 (y) = 0$ for all $y$.
The asymptotic expansion of
the ILDM is given by
Eq.~(\ref{psi-exp}).
A comparison of the coefficients
in the expansion with the coefficients
in the expansion of the slow manifold
$\mathcal{M}_\eps$,
Eq.~(\ref{h-eps-exp}),
shows that
\begin{eqnarray}
  \psi^{(0)} &=& h_0 , \\
  \psi^{(1)} &=& h^{(1)} , \\
  \psi^{(2)} &=& h^{(2)} - (D_z g)^{-2} (D^2 h_0) (f, f) .
\end{eqnarray}
\end{COROLLARY}

The difference $\psi^{(2)} - h^{(2)}$
involves the bilinear form $D^2 h_0$,
which is proportional to the curvature
of the zero set of $g$ at $\eps = 0$.
It is present in any fast-slow system,
unless the curvature vanishes everywhere.
Because of it, the ILDM is in general not invariant
under the dynamics of the system of
Eqs.~(\ref{eq-y})--(\ref{eq-z}).

\section{The Iterative Method of Fraser and Roussel \label{s-FR}}
\setcounter{equation}{0}
The iterative method of
Fraser and Roussel was developed
originally for planar fast-slow systems
that are linear in the fast variable,
\begin{eqnarray}
  \dot{y} &=& f_1 (y, \eps) z  + f_2 (y, \eps), \Label{eq-y-planar} \\
  \eps \dot{z} &=& g_1 (y, \eps) z + g_2 (y, \eps) \Label{eq-z-planar} .
\end{eqnarray}
Here, $y$ and $z$ are scalar-valued functions of time.
These systems of equations are typical
for enzyme kinetics~\cite{F-1975} and
other biochemical systems
whose dynamics can be reduced to
slow manifolds.
In this case, the
slow manifolds are curves
in the phase plane.

On any trajectory
$z = z (y, \eps)$
in the phase plane,
we have the identity
$\dot{z} = z_y \dot{y}$, or,
in terms of the functions
$f$ and $g$,
\be
  \eps z_y (f_1 z + f_2) = g_1 z + g_2 .
  \Label{inv-planar}
\ee
This identity is known as the
\emph{invariance equation}.
In the present case,
it can be solved for $z$
in terms of $y$ and $z_y$,
\be
  z
  = \frac{- g_2 + \eps f_2 z_y}
         {g_1 - \eps f_1 z_y} .
  \Label{inv-planar-z}
\ee
The equation holds, in particular,
along trajectories on invariant manifolds.
Fraser used Eq.~(\ref{inv-planar-z})
to propose the following functional
iteration procedure
to approximate the slow manifold.

Starting from an initial function $\varphi^{(0)}$,
one computes a sequence of functions
$\{ \varphi^{\ell} : \ell = 1, 2, \ldots \}$
using the definitions
\be
  \varphi^{(\ell)}
  = \frac{-g_2 + \eps f_2 \varphi^{(\ell - 1)}_y}
         {g_1 - \eps f_1 \varphi^{(\ell - 1)}_y} , \quad
  \ell = 1, 2, \ldots\,.
  \Label{FR-planar}
\ee
Under appropriate conditions,
the sequence
$\{ \varphi^{(\ell)} (y, \eps) : \ell = 1, 2, \ldots \}$
approaches $z (y, \eps)$
(in a sense to be made precise)
as $\ell$ goes to infinity,
so the algorithm generates
successive approximations
to a slow manifold.

The iterative procedure
generalizes to the fast-slow system
of Eqs.~(\ref{eq-y})--(\ref{eq-z}).
The invariance equation is
\be
  \eps (D_y z) (y, \eps) f (y, z(y, \eps), \eps)
  = g (y, z(y, \eps), \eps) ,
  \Label{inv-gen}
\ee
for any trajectory
$z = z (y, \eps)$ in phase space.
Starting from a function $\varphi^{(0)}$,
one computes a sequence of functions
$\{\varphi^{(\ell)} : \ell = 1, 2, \ldots \}$
by solving the equation
\be
  \eps (D_y \varphi^{(\ell - 1)}) (y, \eps)
  f (y, \varphi^{(\ell)} (y, \eps), \eps)
  =
  g (y, \varphi^{(\ell)} (y, \eps), \eps) .
  \Label{FR-gen}
\ee
The sequence
$\{ \varphi^{(\ell)} (y, \eps) : \ell = 1, 2, \ldots \}$
approaches $z (y, \eps)$
(again, in a sense to be made precise)
as $\ell$ goes to infinity.

Notice that Eq.~(\ref{FR-gen}) amounts to
an \emph{implicit definition} of
$\varphi^{(\ell)}$, unless both
$f$ and $g$ are linear in the
fast variable~$z$,
as in the planar case discussed above,
Eqs.~(\ref{eq-y-planar})--(\ref{eq-z-planar}).
Hence, the numerical computation
of $\varphi^{(\ell)}$ generally
requires the solution of
a nonlinear equation.

\section{Asymptotics of the Iterative Method \label{s-FR-as}}
\setcounter{equation}{0}
Beacuse the iterative method
of Fraser and Roussel
is closely related to
the invariance equation,
its asymptotic properties are
most easily analyzed in terms
of those of Eq.~(\ref{inv-gen}).

\begin{LEMMA} \Label{l-inv-as}
The invariance equation,
Eq.~(\ref{inv-gen}), admits
an asymptotic solution in the form
of a power series expansion,
\be
   z (y,\eps)
   = z^{(0)} (y) + \eps z^{(1)} (y) + \cdots
  \quad \mbox{as } \eps \downarrow 0 ,
  \Label{z-exp}
\ee
where
\be
  z^{(0)} = h_0 , \Label{z0}
\ee
and the functions
$z^{(l)} : K \to \mathbf{R}^n$,
$l = 1, 2, \ldots\,$,
are found successively
from Eq.~(\ref{zl}) below.
In particular, $z^{(1)}$ and $z^{(2)}$
are found from the equations
\begin{eqnarray}
  (D_z g) z^{(1)}
  &=& (Dz^{(0)}) f - g_\eps ,
             \Label{z1} \\
  (D_z g) z^{(2)}
  &=& (Dz^{(1)}) f
     + (Dz^{(0)}) \left( (D_z f) z^{(1)} + f_\eps \right) \nonumber \\
  &&\mbox{}- \half (D_z^2 g) (z^{(1)}, z^{(1)})
                  - (D_z g_\eps) z^{(1)}
                  - \half g_{\eps\eps} ,
             \Label{z2}
\end{eqnarray}
where $f$ and $g$
and their derivatives
are evaluated at
$(y, z^{(0)}(y), 0)$.
\end{LEMMA}

\begin{PROOF}
We begin by expanding the function $f$,
\be
  f (y, z(y, \eps), \eps)
  = f^{(0)} (y) + \eps f^{(1)} (y) + \eps^2 f^{(2)} (y) + \cdots \,,
  \Label{f-exp}
\ee
where
\be
  f^{(0)} (y) = f (y, z^{(0)} (y), 0)
  \Label{f0}
\ee
and
\begin{eqnarray}
  f^{(l)} (y)
  &=& \sum_{j=1}^l \frac{1}{j!} (D_z^j) f (y, z^{(0)} (y), 0)
  \sum_{|i| = l} (z^{(i_1)} (y), \ldots\,, z^{(i_j)} (y)) \nonumber \\ 
  &&\mbox{}+ \sum_{k=1}^{l-1} \frac{1}{k!} 
             \sum_{j=1}^{l-k} \frac{1}{j!} (D_z^j (\partial_\eps^k f)) (y, z^{(0)} (y), 0)
  \sum_{|i| = l-k} (z^{(i_1)} (y), \ldots\,, z^{(i_j)} (y)) \nonumber \\ 
  &&\mbox{}+ \frac{1}{l!} (\partial_\eps^l f) (y, z^{(0)}(y), 0),
   \quad
  l = 1, 2, \ldots\,.
  \Label{fl}
\end{eqnarray}
The derivative $(D_z^j f) (y, z^{(0)} (y), 0)$ in the first term
is a $j$-linear map from
$(\mathbf{R}^n)^j$ to $\mathbf{R}^m$.
The inner sums are taken over all multiindices
$i = (i_1, \ldots\,, i_j)$ 
of positive integers $i_1, \ldots\,, i_j$
with length $|i| = \sum_{k=1}^{j} i_k = l$ 
and $l-k$, respectively.
The first few coefficients are
\begin{eqnarray}
  f^{(1)} (y) &\hspace{-0.5em}=\hspace{-0.5em}& (D_z f) z^{(1)} 
                + f_\eps , \Label{f1} \\
  f^{(2)} (y) &\hspace{-0.5em}=\hspace{-0.5em}& (D_z f) z^{(2)} 
                + \half (D_z^2 f) (z^{(1)}, z^{(1)})
                + (D_z f_\eps) z^{(1)} 
                + \half f_{\eps\eps}  , \Label{f2} \\
  f^{(3)} (y) &\hspace{-0.5em}=\hspace{-0.5em}& (D_z f) z^{(3)} 
                + (D_z^2 f) (z^{(1)}, z^{(2)})
                + \onesixth (D_z^3 f) (z^{(1)}, z^{(1)}, z^{(1)}) \nonumber \\
                &&\mbox{}
                + (D_z f_\eps ) z^{(2)} 
                + \half (D_z^2 f_\eps ) (z^{(1)}, z^{(1)}) \nonumber \\
                &&\mbox{}
                + \half (D_z f_{\eps\eps} ) z^{(1)}
                + \onesixth f_{\eps\eps\eps}  , \Label{f3}
\end{eqnarray}
where $f$ and its derivatives
are evaluated at
$(y,z^{(0)}(y),0)$,
and the argument of each 
$z^{(i)}, i=1,2,\ldots\, ,$ 
is $y$.
A similar expansion holds for $g (y, z(y, \eps), \eps)$,
\be
  g (y, z(y, \eps), \eps)
  = g^{(0)} (y) + \eps g^{(1)} (y) + \eps^2 g^{(2)} (y) + \cdots \,,
  \Label{g-exp}
\ee
where
\be
  g^{(0)} (y) = g (y, z^{(0)} (y), 0) ,
  \Label{g0}
\ee
and
\begin{eqnarray}
  g^{(l)} (y)
  &=&
  \sum_{j=1}^l \frac{1}{j!} (D_z^j g) (y, z^{(0)} (y), 0)
  \sum_{|i| = l} (z^{(i_1)} (y), \ldots\,, z^{(i_j)} (y)) \nonumber \\  
  &&\mbox{}+ \sum_{k=1}^{l-1} \frac{1}{k!} 
             \sum_{j=1}^{l-k} \frac{1}{j!} (D_z^j (\partial_\eps^k g)) (y, z^{(0)} (y), 0)
  \sum_{|i| = l-k} (z^{(i_1)} (y), \ldots\,, z^{(i_j)} (y)) \nonumber \\
  &&\mbox{}+ \frac{1}{l!} (\partial_\eps^l g) (y, z^{(0)}(y), 0),
  \quad l = 1, 2, \ldots\,.
  \Label{gl}
\end{eqnarray}
Termwise differentiation of the
asymptotic expansion~(\ref{z-exp})
gives
\be
  (D_y z) (y, \eps)
  = Dz^{(0)} + \eps Dz^{(1)} + \cdots \,.
  \Label{dy-z-exp}
\ee
Equating the coefficients
of like powers of $\eps$
in the left and right members 
of the invariance equation, Eq.~(\ref{inv-gen}),
we obtain a sequence of functional identities,
\begin{eqnarray}
  g^{(0)} &=& 0 , \Label{eq-g0} \\
  g^{(\ell)} &=& \sum_{i=0}^{\ell-1} (Dz^{(i)}) f^{(\ell-1-i)} , \quad
  \ell = 1, 2, \ldots\,. 
  \Label{eq-gl}
\end{eqnarray}
We satisfy the ${\cal O}(1)$ equation,
Eq.~(\ref{eq-g0}), by taking
$z^{(0)} = h_0$;
see Eq.~(\ref{z0}).
Then
$f^{(0)} (y) = f (y, h_0(y), 0)$, and
$f^{(\ell)}(y)$ and 
$g^{(\ell)}(y)$ 
are given by Eqs.~(\ref{fl})
and (\ref{gl}), respectively,
with $z^{(0)}(y)$ 
replaced by $h_0(y)$.

Next,
we turn 
to the ${\cal O}(\eps^\ell)$ equation,
Eq.~(\ref{eq-gl}).
We observe that 
$z^{(\ell)}$ occurs 
in $g^{(\ell)}$
only in the first term with $j=1$; 
the remaining terms 
involve $z^{(1)}$ 
through $z^{(\ell-1)}$
but not $z^{(\ell)}$.
The right member
of Eq.~(\ref{eq-gl}) 
similarly involves $z^{(1)}$ 
through $z^{(\ell-1)}$
but not $z^{(\ell)}$.
Therefore,
the identities~(\ref{eq-gl}) 
can be solved successively
for $z^{(1)}$, $z^{(2)}$, 
and so on.
Thus we find
\begin{eqnarray}
  (D_z g) z^{(\ell)}
  &\hspace{-0.5em}=\hspace{-0.5em}&
      \sum_{i=0}^{\ell-1} (Dz^{(i)}) f^{(\ell-1-i)} 
       - \sum_{j=2}^\ell \frac{1}{j!} (D_z^j g)
         \sum_{|i| = \ell} (z^{(i_1)}, \ldots \,, z^{(i_j)})  
      \nonumber \\
   &&\hspace{-2em}\mbox{}-
         \sum_{k=1}^{\ell-1} \frac{1}{k!} 
         \sum_{j=1}^{\ell-k} \frac{1}{j!} (D_z^j (\partial_\eps^k g))
                 \sum_{|i| = \ell-k} (z^{(i_1)}, \ldots\,, z^{(i_j)})
         - \frac{1}{\ell!} (\partial_\eps^\ell g)
  \Label{zl}
\end{eqnarray}
for $\ell = 1, 2, \ldots\,$.
Here, the functions 
$f$ and $g$ 
and their derivatives
are evaluated at $(y, z=h_0(y), 0)$,
and it is understood 
that a sum is empty 
when the lower bound 
exceeds the upper bound.
The equations for $\ell=1$ and $\ell=2$
are given in the statement of the theorem.
\end{PROOF}

The following theorem shows that
$\ell$~successive applications of
the iterative algorithm of Fraser and Roussel,
starting from $\varphi^{(0)} = h_0$,
generate an approximation $\varphi^{(\ell)}$
to the slow manifold $\mathcal{M}_\eps$
that is accurate up to and including
the $\mathcal{O}(\eps^\ell)$ term
as $\eps \downarrow 0$.

\begin{THEOREM} \Label{t-FR}
Let $\varphi^{(\ell)}$ and $z^{(\ell)}$
be defined recursively for $\ell = 1, 2,\ldots$
by Eq.~(\ref{FR-gen}) and~(\ref{zl}),
respectively.
If $\varphi^{(0)} = z^{(0)} = h_0$, then
\be
  \varphi^{(\ell)}
  \equiv
  \varphi^{(\ell)} (y, \eps)
  = \sum_{i=0}^\ell \eps^i z^{(i)} (y) + \mathcal{O}(\eps^{\ell+1}) ,
  \quad \ell = 1, 2, \ldots\,.
\ee
\end{THEOREM}

\begin{PROOF}
The proof is by induction.
Taking $\ell=1$, we have
\[
  \eps (D \varphi^{(0)}) (y) f (y, \varphi^{(1)} (y, \eps), \eps)
  = g (y, \varphi^{(1)} (y, \eps), \eps) .
\]
Since $\varphi^{(0)} = z^{(0)}$,
this equation is the same as
\[
  \eps (D z^{(0)}) (y) f (y, \varphi^{(1)} (y, \eps), \eps)
  = g (y, \varphi^{(1)} (y, \eps), \eps) .
\]
We expand the terms in this equation
in powers of $\eps$ and equate
the coefficients of like powers
of $\eps$.
To leading order, we find the equation
\[
  g (y, \varphi^{(1)} (y, 0), 0) = 0 ,
\]
which is precisely the equation for $z^{(0)}$,
so
\[
  \varphi^{(1)} (y, 0) = z^{(0)} (y) .
\]
To the next order, we find the equation
\[
  (D_z g) (y, z^{(0)}, 0) \varphi^{(1)}_\eps (y, 0)
  + g_\eps (y, z^{(0)}, 0)
  = (D z^{(0)}) (y) f (y, z^{(0)} (y), 0) ,
\]
which is precisely Eq.~(\ref{z1}),
so
\[
  \varphi^{(1)}_\eps (y, 0) = z^{(1)} (y) .
\]
Thus,
\[
  \varphi^{(1)} (y, \eps)
  = z^{(0)} (y)
  + \eps z^{(1)} (y) + \mathcal{O}(\eps^2) ,
\]
and the theorem is true for $\ell = 1$.

Suppose the theorem is true for $\ell - 1$,
$\varphi^{(\ell-1)}
= \sum_{i=0}^{\ell-1} \eps^i z^{(i)}
+ \mathcal{O}(\eps^\ell)$.
The function $\varphi^{(\ell)}$
is defined by Eq.~(\ref{FR-gen}),
\[
  \eps (D_y \varphi^{(\ell-1)}) (y, \eps)
  f (y, \varphi^{(\ell)} (y, \eps), \eps)
  = g (y, \varphi^{(\ell)} (y, \eps), \eps) .
\]
We expand each term in powers of $\eps$
and equate the coefficients of like powers,
\[
  \sum_{i=0}^{j-1} (D_y z^{(i)}) f^{(j-1-i)}
  =
  g^{(j)} ,
  \quad j = 1, 2, \ldots\,, \ell ,
\]
where
$f^{(\cdot)}$ and $g^{(\cdot)}$
are the functions
defined after Eq.~(\ref{f-exp}).
For $j=1, \ldots, \ell$,
these are exactly the same functional identities
as we derived above.
Hence, order by order,
the solution is
\[
  \partial_\eps^i \varphi^{(\ell)} (y, 0)
  = z^{(i)} (y) , \quad
  i = 0, 1, \ldots, \ell .
\]
This proves that the theorem is true for $\ell$.
\end{PROOF}

\noindent\textbf{Remark.}
In general, the coefficient of $\eps^2$ in
$\varphi^{(1)}$ will not be equal to $z^{(2)}$.
This may be seen by direct examination of
the $\mathcal{O}(\eps^2)$ equation,
\[
  (D z^{(0)}) (y) \left( (D_z f) z^{(1)} + f_\eps \right)
  =
  \half (D_z^2 g) (z^{(1)}, z^{(1)})
  + (D_z g_\eps) z^{(1)} + \half g_{\eps\eps} ,
\]
which differs from Eq.~(\ref{z2}).
The latter has two additional terms,
namely, $(D_z g) z^{(2)}$ and
$(D z^{(1)}) f$.
Hence, the first iterate generally involves
an error of $\mathcal{O}(\eps^2)$.
Similarly, 
the $\ell$th iterate
does not give the same equation
for the term at order $\mathcal{O}(\eps^{\ell+1})$ 
as compared with that obtained from invariant manifold theory;
hence, the error in the approximation at this stage is 
generally $\mathcal{O}(\eps^{\ell+1})$.~\qed

A comparison of the results
given in Theorem~\ref{t-FR}
with Theorem~\ref{t-Fenichel}
leads to the following conclusions.

\begin{COROLLARY}  \Label{c-FR}
The iterative method of Fraser and Roussel
gives successively higher-order asymptotic
approximations to the slow
manifold~$\mathcal{M}_\eps$.
Starting from $\varphi^{(0)} = h_0$,
$\ell$ applications of the iterative
procedure give an approximation
$\varphi^{(\ell)}$ that satisfies
\be
  \varphi^{(\ell)}
  = \sum_{i=0}^\ell \eps^i h^{(i)}
  + \mathcal{O} (\eps^{\ell + 1}) , \quad
  \ell = 1, 2, \ldots\,,
\ee
where the functions $h^{(i)}$ are the
coefficients in the asymptotic expansion
of $\mathcal{M}_\eps$ given in
Theorem~\ref{t-Fenichel}, 
Eq.~(\ref{h-eps-exp}).
\end{COROLLARY}

\noindent\textbf{Remark.}
In Ref.~\cite{RF-1991a},
Roussel and Fraser
noted the
connection between the iterative approach
developed by Fraser and Roussel and techniques
from dynamical systems theory, such as deriving
the existence of certain invariant manifolds
via application of the contraction mapping
principle to appropriate integral equations.~\qed

\section{The Michaelis--Menten--Henri Model \label{s-MMH}}
\setcounter{equation}{0}
In this section,
we illustrate the analytical results
of the preceding sections on
the Michaelis--Menten--Henri (MMH)
model.
The MMH model is a prototype 
reaction mechanism for enzyme kinetics
in biochemistry~\cite{HTA-1967}.
It is a planar fast-slow system,
which is given in nondimensional form
by the equations
\begin{eqnarray}
  \dot y &=& - y + (y + a - b) z ,
  \Label{MMH-y} \\
  \eps \dot z &=& y - (y + a) z .
  \Label{MMH-z}
\end{eqnarray}
The variables are
$y$, the concentration
of the substrate,
and $z$, the concentration of
an intermediate substrate-enzyme complex.
The parameters $a$ and $b$
satisfy the inequalities $a > b > 0$.

\noindent\textbf{Remark.}
The MMH model,
Eqs.~(\ref{MMH-y})--(\ref{MMH-z}),
is derived from a more complicated
system involving four species
(enzyme, substrate,
enzyme-substrate complex,
and product) and two reactions
(one reversible, one irreversible).
The full system can be reduced to
the planar system, because it has
two conserved quantities;
see Ref.~\cite{NF-1989}.~\qed

The system of 
Eqs.~(\ref{MMH-y})--(\ref{MMH-z})
has a family 
of slow manifolds 
$\mathcal{M}_\eps$,
whose asymptotics 
are given by Eq.~(\ref{h-eps-exp}),
\be
  h_\eps(y)
  = h_0 (y) 
  + \eps h^{(1)} (y)
  + \eps^2 h^{(2)} (y) + \cdots \,,
  \quad y > 0 ,
  \Label{MMH-M-eps}
\ee
where 
\begin{eqnarray}
  h_0 (y) &=& \frac{y}{y+a},
  \Label{h0-MMH} \\
  h^{(1)} (y) &=& \frac{aby}{(y + a)^4} ,
  \Label{h1-MMH} \\
  h^{(2)} (y) &=& \displaystyle
    {\frac{aby (2ab - 3by - ay - a^2)}{(y + a)^7} } .
  \Label{h2-MMH}
\end{eqnarray}
We now show that
the ILDM method 
yields an approximation
of the slow manifold
that is accurate up to and including
the $\mathcal{O}(\eps)$ term,
but differs at $\mathcal{O} (\eps^2)$
because of the curvature of $h_0$.

The Jacobian of the vector field
associated with
Eqs.~(\ref{MMH-y})--(\ref{MMH-z}) is
\be
  J
  = \left(
  \begin{array}{cc}
  - (1 - z)         & y + a - b \\
  \eps^{-1} (1 - z) & - \eps^{-1} (y + a)
  \end{array}
  \right) ,
  \Label{jac-MMH}
\ee
and its eigenvalues are
\be
  \lambda_{s,f} (y, z)
  =\mbox{} - \frac{1}{2} \left(\frac{y + a}{\eps} + 1 - z \right)
  \pm \sqrt{ \frac{1}{4} \left(\frac{y + a}{\eps} + 1 - z \right)^2
  - \frac{b(1 - z)}{\eps}} .
  \Label{lambdapm-MMH}
\ee
The (nonnormalized) slow eigenvector is
\be
  v_s
  = \left(
  \begin{array}{c}
     \lambda_s + \eps^{-1} (y+a) \\
     \eps^{-1} (1-z)
  \end{array}
  \right) ,
  \Label{vs-MMH}
\ee
and there is a corresponding fast eigenvector.
The vector $v_s$ spans the slow subspace.
As noted before, the fast and slow subspace
are not orthogonal, so we work with $v_s^\perp$
and define the ILDM as the set of points
where the vector field is orthogonal to $v_s^\perp$,
\be
  (1-z) [ -y + (y + a - b) z ]
  - (\lambda_{s} + \eps^{-1} (y + a))[y - (y+a) z] = 0 .
  \Label{MPE-MMH}
\ee
Asymptotically, the ILDM
is given by
\be
  z = \psi (y,\eps)
  = \psi^{(0)}(y) + \eps \psi^{(1)}(y) + \eps^2 \psi^{(2)}(y) + \cdots\,,
  \Label{MMH-ILDM}
\ee  
where
\begin{eqnarray}
  \psi^{(0)} (y) &=& \frac{y}{y + a} , \Label{MP-z0-MMH} \\
  \psi^{(1)} (y) &=& \frac{aby}{(y + a)^4} , \Label{MP-z1-MMH} \\
  \psi^{(2)} (y) &=& \frac{aby(2ab - by - ay - a^2)}
                          {(y + a)^7} . \Label{MP-z2-MMH}
\end{eqnarray}
A comparison of the coefficients
in the expansions~(\ref{MMH-M-eps})
and~(\ref{MMH-ILDM}) shows
agreement of the $\mathcal{O}(1)$
and $\mathcal{O}(\eps)$ terms.
On the other hand,
the $\mathcal{O}(\eps^2)$ terms
differ;
their difference is
proportional to the
curvature of $h_0$,
\be
  \psi^{(2)} - h^{(2)}
  = \frac{2ab^2y^2}{(y+a)^7}
  = -\frac{f^2}{g_z^2} h_0'' .
\ee
The iterative method of
Fraser and Roussel starts from
the invariance equation,
\be
  z = \frac{y + \eps z_y y}{y + a + \eps z_y (y + a - b)}
\ee
or, equivalently,
\begin{eqnarray}
  z &=& \frac{y}{y + a}
  + \eps \frac{b y z_y}{(y + a)^2}
  - \eps^2 \frac{b y (y + a - b) z_y^2}{(y + a)^3} \nonumber \\
  &&\mbox{}+
  \eps^3 \frac{b y (y + a - b)^2 z_y^3}
              {(y + a)^4 + \eps (y + a)^3 (y + a - b) z_y} .
\end{eqnarray}
Successive applications of the iterative algorithm
lead to the approximations
\begin{eqnarray}
  \varphi^{(0)} &\hspace{-0.5em}=\hspace{-0.5em}& \frac{y}{y+a} , \\
  \varphi^{(1)} &\hspace{-0.5em}=\hspace{-0.5em}& \frac{y}{y+a} + \eps \frac{aby}{(y+a)^4}
  - \eps^2 \frac{a^2 b y (y+a-b)}{(y+a)^7} + \mathcal{O} (\eps^3) , \\
  \varphi^{(2)} &\hspace{-0.5em}=\hspace{-0.5em}& \frac{y}{y+a} + \eps \frac{aby}{(y+a)^4}
  + \eps^2 \frac{aby (2ab - 3by - ay - a^2)}
              {(y+a)^7} + \mathcal{O} (\eps^3) .
\end{eqnarray}
A comparison with Theorem~\ref{t-Fenichel}
(and Eqs.~({\ref{h0-MMH})--(\ref{h2-MMH}))
shows that $\varphi^{(\ell)}$
is asymptotically correct
up to and including terms of
$\mathcal{O} (\eps^\ell)$
for $\ell = 1, 2$ (and beyond),
as predicted by the analysis.

\noindent\textbf{Remark.}
The analysis of the full MMH model
has been significantly extended
in Ref.~\cite{SS-1989}.~\qed

\section{The Davis-Skodje Model \label{s-DS}}
\setcounter{equation}{0}
The planar fast-slow system 
\begin{eqnarray}
   \dot y & = & -y , \Label{ds-eq-y} \\
   \eps \dot z & = & -z + \frac{y}{1+y} - \frac{\eps y}{(1+y)^2},
   \Label{ds-eq-z}
\end{eqnarray}
was introduced by Davis and Skodje~\cite{DS-1999}
as a model on which to compare 
various reduction methods. 
(The inverse $\eps^{-1}$,
which is large,
equals the large parameter 
$\gamma$ of Ref.~\cite{DS-1999}.)

For any $\eps$,
the curve
\be
  z = h_\eps (y) = \frac{y}{1+y} , \quad y \geq 0 ,
  \Label{slowmfld-ds-eqn}
\ee
is invariant 
under the dynamics 
of Eqs.~(\ref{ds-eq-y})--(\ref{ds-eq-z}).
Therefore, the function $h_\eps$
represents the slow manifold
exactly on $y \geq 0$
for all small $\eps > 0$.
(The nonlinearity in
Eqs.~(\ref{ds-eq-y})--(\ref{ds-eq-z})
was, in fact, chosen so
the slow manifold is given
by the simple expression
of Eq.~(\ref{slowmfld-ds-eqn}).)

The Jacobian of the vector field is
\be
  J = \left(
      \begin{array}{cc}
         -1 & 0 \\
         \eps^{-1} ((1+y) + \eps (y-1)) (1+y)^{-3} & \mbox{}- \eps^{-1}
      \end{array}
      \right) ,
  \label{Jac-ds-eqn}
\ee
and the eigenvalues are
$\lambda_s = -1$ and
$\lambda_f = -1/\eps$
for all $(y, z)$.
The corresponding
(nonnormalized) eigenvectors are
\be
  v_s
  = \left( 
         \begin{array}{c}
             1 \\
             (1 - \eps)^{-1} (1 + y + \eps (y-1)) (1 + y)^{-3}
         \end{array}
         \right) , \quad
  v_f
  = \left(
         \begin{array}{c}
             0 \\
             1
         \end{array}
         \right) .
\ee
The ILDM is given by
the expression
\be
    z = \psi (y, \eps)
    = \frac{y}{1+y}
    + \frac{2 \eps^2 y^2}{(1-\eps)(1+y)^3} ;
    \Label{ILDM-ds-eqn}
\ee
cf.~\cite[Eq.~(3.8)]{DS-1999}.
Its asymptotic expansion is
\be
  z
  =
  \frac{y}{1+y}
  + \eps^2 \frac{2y^2}{(1+y)^3}
  + \cdots \,.
  \Label{asymp-ILDM-ds-eqn}
\ee
The error in the expansion
is $\mathcal{O}(\eps^2)$
and proportional to
the curvature,
$h_0'' = -2/(1+y)^3$.

The invariance equation is
\be
  z = \frac{y}{1+y} - \eps \frac{y}{(1+y)^2} + \eps y z_y ,
\ee
from which one readily verifies that
the iterative method of Fraser and Roussel
yields the approximation
$\varphi^{(\ell)} (y) = y/(1+y)$
for $\ell = 1, 2, \ldots\,$.

\noindent\textbf{Remark.}
If one restricts 
the variable $y$
to a finite interval,
then the system of
Eqs.~(\ref{ds-eq-y})--(\ref{ds-eq-z})
has a family of slow manifolds,
all exponentially close 
($\mathcal{O}(\mathrm{e}^{-c/\eps)})$
for some $c>0$)
to the exact slow manifold,
Eq.~(\ref{slowmfld-ds-eqn}).~\qed

\section{Discussion} \label{s-disc}
\setcounter{equation}{0}
The fast-slow system of
Eqs.~(\ref{eq-y})--(\ref{eq-z})
captures the essential elements
of any reaction mechanism whose
long-time dynamics evolve on
a slow manifold in the composition space.
As stated in Section~\ref{s-general}, however,
this system is a mathematical idealization,
and we need to consider how the results
of the analysis carry over to
more general reaction mechanisms.
In this subsection,
we consider issues related to
the separation of time scales
(Section~\ref{s-disc-timescales}),
the inclusion of conserved quantities
(Section~\ref{s-disc-conserved}),
and the development and analysis 
of reduction mechanisms for
reaction-diffusion equations
(Section~\ref{s-disc-diffusion}).

\subsection{Separation of Time Scales \label{s-disc-timescales}}
In this section,
we discuss the partition
of variables into
a fast and a slow group
and the assumption that the groups
evolve on time scales that are
and remain well separated at all times.
This assumption underlies the definition
of the small parameter~$\eps$ in
the model of Eqs.~(\ref{eq-y})--(\ref{eq-z}).
We also discuss the possibility
of partitioning the variables
in more than two groups.

While many of the systems
to which reduction methods have been applied
satisfy this assumption,
there are a significant number
of reaction mechanisms
where the fast and slow
time scales are separated,
but not well separated.
This is the case, for example,
when $\eps$ is no longer an
asymptotically small parameter
but a fixed (relatively small) number.
In such cases,
the spectral gap between
the fast and slow eigenvalues
of the Jacobian of the vector field
is small,
much smaller than the chasm between
the fast $\mathcal{O}(\eps^{-1})$ eigenvalues
and the slow $\mathcal{O}(1)$ eigenvalues
for asymptotically small $\eps$,
even at points near a low-dimensional manifold.
Nevertheless,
as long as there exists a nonzero spectral gap,
the Jacobian can still be reduced 
to a fast and a slow component,
and the ILDM method can be 
(and has been) implemented numerically.
Part of our ongoing research
is aimed at using spectral projection operators
to analyze these applications of the ILDM method.

In addition,
there are reaction mechanisms
where the number of fast and
slow species changes over time,
as may happen, for example,
when the temperature in
the chemical reactor changes
and the least-slow 
of the slow species 
transits to the group
of fast species.
To account for this type
of occurrence,
practical implementations
of the ILDM method may use
one ILDM until
the crossover occurs
and another after the crossover.
The two procedures can be linked
by doing a numerical integration
of the full system of equations
during the crossover.
The analysis presented in
this article applies
to each ILDM separately.
However,
it should be noted that
the relationship between 
these two ILDMs
(and, more generally, slow manifolds)
depends on the bifurcation
that occurs at crossover.
As an alternative,
one may consider repartitioning
the species to avoid crossover altogether.

In some systems,
a component of $y$
may evolve on an even slower time scale
than that given by $\tau=\eps t$.
For example, in Eq.~(\ref{eq-y})
one may have
$y_i' = \eps^{1 + \gamma} f_i (y, z, \eps)$
with $\gamma > 0$ for some index $i$.
Such cases are accounted for
by the present analysis;
in fact, it suffices to absorb
the factor $\eps^\gamma$
in $f_i$.
On the other hand,
the fact that there are slower
time scales in the model may
point to the existence of still
lower-dimensional slow manifolds.
By systematically eliminating
fast variables, starting with the
fastest and proceeding up the
hierarchy, one can reduce the
dimensionality of the slow manifold
in a systematic way until no
further reduction is possible.
Such an approach has been used,
for example, in Refs.~\cite{RF-1991a, YTBMP-1995}.
The idea of a hierarchy of time
scales was first explored
in singular perturbation theory
by Tikhonov~\cite{T-1948}.

\subsection{Conserved Quantities and
  the ILDM Method \label{s-disc-conserved}}
In this section,
we briefly consider
the ILDM method
for fast-slow systems
whose dynamics are described by
Eqs.~(\ref{eq-y})--(\ref{eq-z}),
where the unknowns satisfy one
or more conservation laws.
For models of chemical reactions,
a conserved quantity is,
typically, a linear combination
of several unknowns that
is constant in time.
Conserved quantities
give rise to 
zero eigenvalues 
of the Jacobian of the vector field,
and the ILDM method groups
these zero eigenvalues
with the slow ones.

Degeneracies in the Jacobian
affect the analysis 
of the ILDM method.
Consider, for example,
a system given by
Eqs.~(\ref{eq-y})--(\ref{eq-z})
with $m=2$ and $n=2$,
which has one conserved quantity.
If the conserved quantity
is a linear combination
of the two slow variables,
then the first and second row
of the Jacobian 
are linearly dependent.
In this case,
there is effectively
only one slow variable,
and the analysis 
presented in this article
applies to the one-dimensional 
slow manifold $\mathcal{M}_0$.
If, by contrast,
the conserved quantity
is a linear combination
of the two fast variables,
then not only is $J$ degenerate,
but also the columns of
$D_z g(y, h_0(y), 0)$
are linearly dependent;
hence, $\mathcal{M}_0$
is not asymptotically stable
in the four-dimensional 
composition space.
A reduction of the number
of fast variables by one
lifts this degeneracy,
because $\mathcal{M}_0$
is asymptotically stable
in the three-dimensional
reduced composition space.
Other possibilities are that
a conserved quantity depends
on a mix of fast and slow variables
or that there are multiple
conserved quantities.
In each case, the analysis of the
ILDM method presented here applies
once the system has been reduced
to a system for which the manifold
$\mathcal{M}_0$ is asymptotically
stable.

\subsection{Reaction-Diffusion Equations \label{s-disc-diffusion}}
Reduction methods have also been developed
for systems of \emph{reaction-diffusion
equations}~\cite{L-1993, LG-1994, YTBMP-1995,
YTBMP-1996, HG-1999, Kan-1999, SRPP-2001, LJL-2001}.
The elimination of the fast species
affects not only the reaction kinetics
of the slow species (as is the case 
for the kinetics reduction methods
considered in this article)
but also their diffusivities.
The ``effective'' diffusivities
differ from the ordinary diffusivities
by concentration-dependent terms
that are higher order in $\eps$.

We illustrate this phenomenon on the
Michaelis--Menten--Henri mechanism
with diffusion of the slow species,
\begin{eqnarray}
  \dot{y}
  &=& 
  - y + (y + a - b) z 
  + D \Delta y , \\
  \eps \dot{z}
  &=&
  y - (y+a) z .
  \Label{disc-MMH-PDE}
\end{eqnarray}
The variables $y$ and $z$
depend not only on time
but also on space;
$\Delta$ is the Laplace operator.
Ideas from inertial-manifold theory
have been applied to this
reaction-diffusion system~\cite{YTBMP-1995}.
The slow manifold is infinite dimensional,
its asymptotics are given by
\be
  z
  = \frac{y}{y+a} 
  + \eps \left[ \frac{aby}{(y+a)^4}
              - \frac{a}{(y+a)^3} D \Delta y
           \right]
  + \mathcal{O} (\eps^2) ,
  \Label{disc-MMH-slow-mfld}
\ee
and the reduced reaction-diffusion equation 
up to and including terms of $\mathcal{O}(\eps)$
is
\be
  \dot{y}
  = - y
  + (y+a-b) \left( \frac{y}{y+a} + \eps \frac{aby}{(y+a)^4} \right)
  + D \left( 1 - \eps \frac{a (y+a-b)}{(y+a)^3} \right) \Delta y .
  \Label{disc-MMH-reduced-RD}
\ee
The regular diffusivity $D$ is seen to be
corrected by an $\eps$-dependent term
that involves the concentration
of the slow species.

Singular perturbation theory
provides an alternative method
to find the infinite-dimensional
slow manifold and, hence, the reduced
reaction-diffusion equation for this
and similar systems.
In particular, one finds
an asymptotic expansion
for the slow manifold of the form 
\be
  z
  = h_0(y) + \eps h^{(1)} (y, \Delta y)
  + \eps^2 h^{(2)} (y, \Delta y, \Delta^2 y) + \cdots \, .
  \Label{disc-SPT-PDE}
\ee
For the MMH model,
this procedure leads to
the same expansion
of the slow manifold,
Eq.~(\ref{disc-MMH-slow-mfld}),
and dynamical systems theory
states that $\mathcal{M}_\eps$
is the infinite-dimensional weak stable manifold
of the spatially homogeneous state $(0,0)$.

We remark that one can 
include diffusion
of the fast variable $z$,
but only if the diffusion coefficient
is $\mathcal{O}(\eps)$,
\begin{eqnarray}
  \dot{y} &=& f (y, z, \eps) + D_1 \Delta y ,  \Label{eq-D1-y} \\
  \eps \dot{z} &=& g (y, z, \eps) + \eps D_2 \Delta z . \Label{eq-D2-z}
\end{eqnarray}
In this case,
$z^{(0)}$ is independent of $D_2$,
and one can follow the same
asymptotic procedure as before;
the diffusion coefficient $D_2$
enters into the equation of order
$\mathcal{O}(\eps)$.
However, if the diffusion term
in the fast equation is $\mathcal{O}(1)$,
the asymptotic procedure for
finding $z^{(0)}$ breaks down.

\begin{center}
\textbf{ACKNOWLEDGMENTS}
\end{center}
\noindent
We thank
Simon Fraser, Ulrich Maas,
Marc Roussel, and Rex Skodje 
for generously sharing
their insights into
reduction methods.
We also thank our colleagues
Michael Davis and Paul Fischer
(Argonne National Laboratory)
for many stimulating conversations
in the course of this investigation.

\noindent
The work of H.~K.\ was supported by the Mathematical, Information, and
Computational Sciences Division subprogram of the Office of Advanced
Scientific Computing Research, U.S.~Department of Energy, under
Contract W-31-109-Eng-38.
The work of T.~K.\ was supported in part 
by the Division of Mathematical Sciences
of the National Science Foundation via grant NSF-0072596.
T.~K. also thanks the Mathematics and Computer Science Division
at Argonne National Laboratory for its hospitality
and support.

\noindent
Corresponding author:

\noindent
Hans G.\ Kaper \\
Mathematics and Computer Science Division \\
Argonne National Laboratory \\
Argonne, IL 60439

\noindent
Authors' e-mail addresses:

\noindent
\texttt{kaper@mcs.anl.gov} \\
\texttt{tasso@math.bu.edu} \\

\end{document}